%% file: collocation-review.tex
\documentclass[11pt]{amsart}

\pdfoutput=1
\usepackage[left=1.25in,right=1.25in,top=1in]{geometry}

\usepackage{amsmath,amsthm,latexsym,amssymb,graphicx,bbm,microtype,mathtools}
\usepackage{verbatim}
\usepackage{mathrsfs}
\usepackage{latex-macros}
\newtheorem{remark}{Remark}

\usepackage[T1]{fontenc}
\def\naive{na\"{\i}ve}

\newcommand{\vertiii}[1]{{\left\vert\kern-0.25ex\left\vert\kern-0.25ex\left\vert #1
    \right\vert\kern-0.25ex\right\vert\kern-0.25ex\right\vert}}

\newtheorem{theorem}{Theorem}

\title{Stochastic collocation on unstructured multivariate meshes}

\author{Akil Narayan}
\thanks{Mathematics Department, University of Massachusetts Dartmouth, 285 Old Westport Road, North Dartmouth, MA 02747}
\author{Tao Zhou}
\thanks{Institute of Computational Mathematics and Scientific/Engineering Computing, AMSS, the Chinese Academy of Sciences, Beijing, China} 
\thanks{T. Zhou is supported by the National Natural Science Foundation of China (Award Nos. 91130003 and 11201461)}

\begin{document}
\maketitle
\begin{abstract}
Collocation has become a standard tool for approximation of parameterized sys- tems in the uncertainty quantification (UQ) community. Techniques for least-squares regularization, compressive sampling recovery, and interpolatory reconstruction are becoming standard tools used in a variety of applications. Selection of a collocation mesh is frequently a challenge, but methods that construct geometrically “unstructured” collocation meshes have shown great potential due to attractive theoretical properties and direct, simple generation and implementation. We investigate properties of these meshes, presenting stability and accuracy results that can be used as guides for generating stochastic collocation grids in multiple dimensions.
\end{abstract}

\pagestyle{myheadings}
\thispagestyle{plain}
\markboth{A.~Narayan AND T.~Zhou}{Collocation on unstructured meshes}

\input{content/introduction}
\input{content/regression}
\input{content/cs}
\input{content/interpolation}
\input{content/conclusion}

\bibliographystyle{plain}
\bibliography{collocation-review}

\end{document}

%% file: content/introduction.tex
\section{Introduction}

The field of uncertainty quantification has enjoyed much attention in recent years as theoreticians and practitioners have tackled problems in the diverse areas of stochastic analysis, exascale applied computing, high-dimensional approximation, and Bayesian learning. The advent of high-performance computing has led to an increasing demand for efficiency and accuracy in predictive capabilities in computational models.

One of the persistent problems in Uncertainty Quantification (UQ) focuses on parameterized approximation to differential systems: Let $u$ be a state variable for a system that is the solution to a physical model
\begin{align}\label{eq:model-system}
  \mathcal{L}(u; t, x, \omega) = 0,
\end{align}
Above $x \in \R^p$ for $p=1,2,3$ is a spatial variable, $t \in \R$ is a temporal variable, and $\omega \in \Omega$ is a probabilistic event that encodes randomness on a complete probability space $(\Omega, \mathcal{F}, \mathcal{P})$. We assume that the model \eqref{eq:model-system} defines a map $\omega \mapsto u(t, x, \omega)$ with $u: \Omega \rightarrow B$ that is well-posed almost surely for some appropriate space $B$ of $(x,t)$-dependent functions.

The operator $\mathcal{L}$ may represent any mathematical model of interest; examples that are popular in modern applied communities are elliptic partial differential equations, systems of time-dependent differential equations, parametric inverse problems, and data-driven optimization; e.g., \cite{babuska_stochastic_2010,pulch_generalised_2012,marzouk_stochastic_2009,poette_uncertainty_2009,breidt_measure-theoretic_2011,najm_uncertainty_2009,debusschere_protein_2003}. The system defined by the operator $\mathcal{L}$ may include boundary value constraints, initial value prescriptions, physical domain variability, or any combination of these \cite{xiu_fast_2009,Xiu,XiuK1,tartakovsky_stochastic_2006}.

The sought system response $u(t, x, \omega)$ is random/stochastic, given by the solution to \eqref{eq:model-system}. 
The stochastic dependence in \eqref{eq:model-system} given by the event $\omega$ is frequently approximated by a $d$-dimensional random variable $Z(\omega)$. In some cases this parameterization of randomness is straightforward: e.g., in a Bayesian framework when ignorance about the true value of a vector of parameters is modeled by treating this parameter set as a random vector in \eqref{eq:model-system}. In contrast, it is common in models for an infinite-dimensional random field to contribute to the stochasticity, and in these cases parameterization is frequently accomplished by some finite-dimensional truncation procedure, e.g., via the Karhunen-Loeve expansion \cite{Ghanem,xiu_numerical_2010}, and this reduces the stochastic dependence in \eqref{eq:model-system} to dependence on a random vector $Z$. In either case, a modeler usually wants to take $d \triangleq \dim Z$ as large as possible to encode more of the random variability in the model.

Under an assumption of model validity, the larger the stochastic truncation dimension $d$, the more accurate the resulting approximation. (Even when model validity is suspect, one can devise metamodeling procedures to capture model form error \cite{kennedy_bayesian_2001}.) Therefore, it is mathematically desirable to take $d$ as large as possible. We rewrite \eqref{eq:model-system} to emphasize dependence on the $d$-dimensional random variable $Z$:
\begin{align}\label{eq:parameterized-model-problem}
  \mathcal{L}(u; t, x, Z) = 0.
\end{align}
In this article we are ultimately interested in approximating $u(x,t,Z)$ or some functional of it, and concentrate on the task of approximating $u$ as a function of $Z$. This is the standard \textit{modus operandi} for non-intrusive methods.

A major challenge for modern uncertainty quantification is the curse of dimensionality. Coined by Richard Bellman \cite{bellman_adaptive_1961}, this refers to the exponentially-increasing computational cost of resolving variability with respect to an increasing number of parameters. The trade-off that one frequently makes is that a large $d$ induces an accurate stochastic truncation, but results in a computationally challenging problem since $u$ depends on a $d$-dimensional parameter.

When $Z$ is high-dimensional, model reduction techniques such as proper orthogonal decomposition methods \cite{deane_lowdimensional_1991,burkardt_pod_2006} or reduced basis methods \cite{patera_reduced_2007,prudhomme_reliable_2002} are useful. In many situations these methods are powerful in their own right, robustly addressing problems in the scientific community; however, many implementations of these approaches are intrusive, meaning that significant rewrite of large legacy codebases is required. The focus of this paper is directed towards a different approach: non-intrusive response construction using multivariate polynomial collocation. ``Non-intrusive" effectively means that existing black-box tools can be used in their current form. In particular we will focus on weighted methods, which are of concern for \textit{stochastic collocation} methods. Stochastic collocation entails polynomial approximation in parameter ($Z$) space using either interpolation or regularized collocation approximation. These approaches have become extremely popular \cite{XiuH,xiu_fast_2009,burkardt_comparison_2009} for their efficiency and effectiveness. In many situations of interest, polynomial approximations converge to the true response exponentially with respect to the polynomial degree.

In stochastic collocation, a polynomial surrogate that predicts variability in parameter space is constructed from point-evaluations of the model response \eqref{eq:parameterized-model-problem} at an ensemble of fixed parameter values $Z_n \in D$; we will call this ensemble of parameter values a \textit{grid} or \textit{mesh}, or a collection of \textit{nodes}. While much work exists on geometrically structured meshes (e.g. tensor-product lattices or sparse grids), we will focus on \textit{unstructured} meshes, which we believe is fertile ground deserving of much attention. Our use of terms `structured' versus `unstructured' refers to visual appearance of a lattice or geometric regularity of the mesh distribution in multivariate space. Obviously use of such a term is a subjective matter, and our goal is not to taxonomically classify collocation methods as structured or unstructured. Instead, the goal of this paper is to highlight some recent collocation strategies that distribute collocation nodes in an apparently unstructured manner; many of these recently developed methods produce approximation meshes that have attractive theoretical and computational properties.

\section{Generalized polynomial Chaos}

The generalized Polynomial Chaos method (gPC) \cite{XiuK2} is essentially the strategy of approximating the $Z$-dependence of $u(x,t,Z)$ from \eqref{eq:parameterized-model-problem} by a $Z$-polynomial. Let $Z$ be a random variable with density function $\rho(z)$, so that $P[Z \in A] = \int_A \rho(z) \dx{z}$ for any Borel set $A$ contained in the domain of $Z$, $A \subset D \subset \R^d$.

Hereafter, we will subsume $t$-dependence into the variable $x$ and write $u(x,Z) \triangleq u(x,t,Z)$. A gPC method proceeds by making the ansatz that the variability of $u$ in the random variable $Z$ is described by a polynomial:
\begin{align}\label{eq:u-ansatz}
  u_N(x, Z) = \sum_{n=1}^N \widehat{c}_n(x) \phi_n(Z),
\end{align}
for a prescribed multivariate polynomial basis $\phi_n$ and unknown coefficient functions $\widehat{c}_n(x)$. The gPC approach is to choose the basis $\phi_n$ to be the family of $L^2_\rho$-orthogonal polynomials,
\begin{align*}
  \E \phi_n(Z) \phi_m(Z) &= \int_{\R^d} \phi_n(z) \phi_m(z) \rho(z) \dx{z} = \delta_{n,m}
\end{align*}
with $\delta_{n,m}$ the Kronecker delta function. We also make the assumption that these polynomials are complete in the corresponding $\rho$-weighted $L^2$ space; this is satisfied if, for example, $\rho$ is continuous and decays at least as fast as $\exp(-\left|z\right|)$ as $|z| \rightarrow \infty$. More intricate conditions can be found in \cite{ernst_convergence_2012}.


In many cases of practical interest, it is natural to assume that the components of $Z$ are mutually independent. In this case, the multivariate functions $\phi_n$ decompose into products of univariate functions. If the components of $Z$ are independent, then $D = \prod_{j=1}^d I_j$ for univariate intervals $I_j \subset \R$ and $\rho(z) = \prod_{j=1}^d \rho_j\left(z^{(j)}\right)$, with $z = \left( z^{(1)}, z^{(2)}, \ldots, z^{(d)}\right)$ the components of $z$. Then the multivariate orthogonal polynomials are products of univariate orthogonal polynomials:
\begin{align*}
  \phi_{\alpha}(z) &= \prod_{j=1}^d \phi^{(j)}_{\alpha_j}\left(z^{(j)}\right), & \int_{I_j} \phi^{(j)}_n(z) \phi^{(j)}_m(z) \rho_j(z) \dx{z} = \delta_{n,m}.
\end{align*}
We have introduced multi-index notation: $\alpha = \left( \alpha_1, \ldots, \alpha_d\right) \in \N_0^d$ is a multi-index, $z^\alpha = \prod_{j=1}^d \left(z^{(j)}\right)^{\alpha_j}$, and $|\alpha| = \sum_{j=1}^d \alpha_j$. Depending on the situation, we will alternate between integer and multi-index notation for the polynomial ansatz:
\begin{align*}
  u_N(x, Z) = \sum_{n=1}^N \widehat{c}_n(x) \phi_n(Z) = \sum_{\alpha \in \Lambda} \widehat{c}_\alpha(x) \phi_\alpha(Z),
\end{align*}
where $\Lambda \subset \N_0^d$ is an index set with cardinality $\left|\Lambda\right| = N$. Any convenient mapping between $1, \ldots, N$ and $\Lambda$ may be used. The choice of $\Lambda$ defines the polynomial approximation space.

Some classical polynomial spaces to which $u_N$ may belong are the tensor-product, total degree, and hyperbolic cross spaces, respectively:
\begin{subequations}\label{eq:polynomial-spaces}
\begin{align}
  \Lambda^P_{d,k} &= \left\{ \alpha \in \N_0^d\, | \, \max_j \alpha_j \leq k\right\}, & P^d_k &= \mathrm{span} \left\{ z^\alpha\, | \, \alpha \in \Lambda^P_{d,k} \right\},\\
  \Lambda^T_{d,k} &= \left\{ \alpha \in \N_0^d\, | \, |\alpha| \leq k \right\}, & T^d_k &= \mathrm{span} \left\{ z^\alpha\, | \, \alpha \in \Lambda^T_{d,k} \right\} \\
    \Lambda^H_{d,k} &= \left\{ \alpha \in \N_0^d\, | \, \prod_j \alpha_j+1 \leq k+1 \right\}, & H^d_k &= \mathrm{span} \left\{ z^\alpha\, | \, \alpha \in \Lambda^H_{d,k} \right\}
\end{align}
\end{subequations}
We have chosen particular definitions of these spaces above, but there are generalizations. E.g., dimensional anisotropy can be used to `bias' the index set toward more important dimensions, or a different $\ell^p$ norm ($0 < p < 1$) can be placed on index space to tailor the hyperbolic cross spaces. The dimensions of $T^d_k$ and $P^d_k$ are
\begin{align}\label{eq:space-dimensions}
  t^d_k &\triangleq \dim T^d_k = \left( \begin{array}{c} d+k \\ k \end{array}\right), &
  p^d_k &\triangleq \dim P^d_k = (k+1)^d.
\end{align}
The dimension of $H^d_k$ has the following upper bound \cite{FabioL2}:
\begin{align*}
h^d_k \doteq \dim H^d_k \leq \lfloor (k+1) (1+\log(k+1))^{d-1} \rfloor.
\end{align*}
For index sets $\Lambda$ that do not fall into the categories defined by \eqref{eq:polynomial-spaces}, we will use the notation $P(\Lambda)$ to denote the corresponding polynomial space.

\begin{figure}
\begin{center}
  \includegraphics[width=\textwidth]{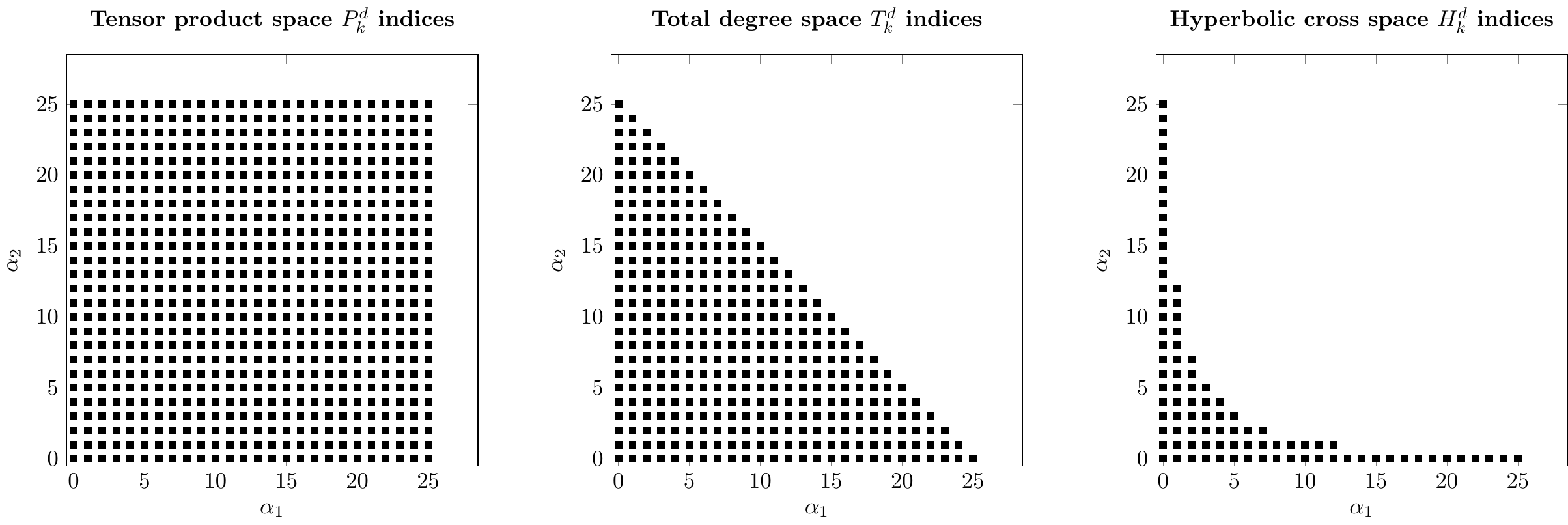}
\end{center}
\caption{Indices associated with the tensor product set $\Lambda^P_{d,k}$ (left), the total degree set $\Lambda^T_{d,k}$ (center), and the hyperbolic cross set $\Lambda^H_{d,k}$ (right). All sets have dimension $d=2$ and degree $k = 25$.}\label{fig:active-indices}
\end{figure}

For the index sets \eqref{eq:polynomial-spaces}, we immediately see the curse of dimensionality: the dimensions of $T^d_k$ and $P^d_k$ increase exponentially with $d$, although $t^d_k$ is smaller than $p^d_k$. The indices in the sets $\Lambda^P_{d,k}$, $\Lambda^T_{d,k}$, and $\Lambda^H_{d,k}$ are graphically plotted in Figure \ref{fig:active-indices} for $d=2$ and polynomial degree $k = 25$.

This highlights a challenge with gPC in high dimensions: the number of degrees of freedom required to resolve highly oscillatory structure grows exponentially with dimension. (Indeed, this is a challenge for any non-adapted multivariate approximation scheme.) In the next section we narrow our focus to collocation schemes.

\subsection{Stochastic Collocation}
The determination of $\widehat{c}_n(x)$ in \eqref{eq:u-ansatz} is the main difficulty, and one way to proceed is to ask that at some predetermined realizations $\left\{ z_m \right\}_{m=1}^M$ of $Z$, the ansatz match the actual response:
\begin{align}\label{eq:collocation-conditions}
  u_N(x, z_m) &\approx u(x, z_m), & m&= 1,\ldots, M,
\end{align}
where the approximation $\approx$ and the relation between $N$ and $M$ are discussed later in this section. (We will allow the ensemble of realizations to be randomly generated in some cases, but we will continue to use lowercase notation $z_m$ to denote specific samples of the random variable $Z$.) This is a stochastic collocation approach and is \textit{non-intrusive}: we need to compute $u(x, z_m)$, but this is accomplished by simply setting $Z = z_m$ in \eqref{eq:parameterized-model-problem} and solving $M$ realizations of the model equation. Therefore, existing deterministic solvers can be utilized. This ability to reuse existing solvers is one of the major strengths of non-intrusive (in particular, collocation) strategies.

In contrast, \textit{intrusive} methods generally require a nontrivial mathematical reformulation of the model \eqref{eq:parameterized-model-problem}, and usually necessitate novel algorithm development. One popular intrusive method for gPC is the stochastic Galerkin approach, where $u_N$ is specified by imposing that some probabilistic moments of the model equation residual vanish. Because these moment equations are coupled and are a novel formulation compared to the original type of equation, existing deterministic solvers of \eqref{eq:parameterized-model-problem} cannot be used. The advantage of intrusive methods compared to non-intrusive methods is that one can usually make more formal mathematical statements about convergence of $u_N$ with an intrusive formulation. For details of intrusive methods, see \cite{Ghanem}. Although intrusive methods are advantageous in many situations, in this paper we only consider non-intrusive collocation approaches.

The focus of this paper discusses how to enforce the collocation conditions \eqref{eq:collocation-conditions}. We investigate three situations:
\begin{itemize}
  \item $M > N$: when there are more constraints than degrees of freedom, we employ regression to attain a solution
  \item $M < N$: with more degrees of freedom than constraints, we may seek sparse solutions and appeal to the theory of compressive sampling
  \item $M = N$: with an equal number of linear constraints and degrees of freedom, we may enforce interpolation
\end{itemize}

Examples of sampling strategies that we consider as `structured' are tensor-product constructions and sparse grid constructions. The former is quickly seen as infeasible for large dimensionality $d$. If we have an $m$-point one-dimensional grid (such as a Clenshaw-Curtis grid), then an isotropic tensorization has $M = m^d$ samples; this dependence on $d$ is usually not computationally acceptable.

Sparse grids are unions of anistropically-tensorized grids, and have proven very effective
\cite{FabioC3,bungartz_sparse_2004,gerstner_numerical_1998,barthelmann_high_2000} at approximating high-dimensional problems.  However, the sparse grids’ adherence to rigid and predictable layouts have the potential weakness of ‘miss- ing’ important features that do not line up with cardinal directions or coordinates in multivariate space, especially earlier in the computation when adaptive methods are seeded with isotropic grids. We believe that unstructured sample designs have the potential to mitigate these shortcomings.


\subsection{Multivariate collocation}\label{sec:methods}

We introduce notation that is used throughout this article. In particular, we reserve the notation $N$ and $M$ to denote the number of terms in the expansion \eqref{eq:u-ansatz} and the number of collocation points in the ensemble \eqref{eq:collocation-conditions}, respectively. We will also use $k$ to denote the polynomial degree of an index set $\Lambda$ when this index set corresponds to one of the choices \eqref{eq:polynomial-spaces}.

The main goal is to compute the coefficient functions $\widehat{c}_n$ from \eqref{eq:u-ansatz}, and this can be decomposed into smaller linear algebra problems. In practice the realization $u(x,z_m)$ (as a solution to \eqref{eq:parameterized-model-problem} with $Z = z_m$) is usually  computed via a spatial discretization as an $R$-dimensional vector $\mathbf{u}(z_m)$. The type of spatial discretization to which $\mathbf{u}(z_m)$ corresponds typically does not influence $Z$-approximation strategies if non-intrusive methods are considered. 

Let $\mathbf{A}$ be an $M \times N$ matrix with entries $(\mathbf{A})_{m,n} = \phi_n(x_m)$ with $\phi_n$ the orthogonal gPC basis. Then the conditions \eqref{eq:collocation-conditions} can be written as
\begin{align*}
  \left( \mathbf{A} \otimes \mathbf{I} \right) \left( \mathbf{C} \right) \approx \mathbf{U},
\end{align*}
where $\mathbf{C}$ is an $R N \times 1$ vector with entries $C_{r,n} = \widehat{c}_{n,r}$ ordered lexicographically in $(r,n)$, $\mathbf{U}$ is an $R M  \times 1$ vector with entries $U_{r,m} = \left(\mathbf{u}(z_m)\right)_r$ ordered lexicographically, and $\mathbf{I}$ is the $R \times R$ identity matrix. Then clearly the above system can be rewritten as the following decoupled series of systems:
\begin{align}\label{eq:decoupled-linear-systems}
  \mathbf{A} \mathbf{c}_r &\approx \mathbf{u}_r, & r &= 1, \ldots, R,
\end{align}
where $\mathbf{c}_r$ has entries $\left( \mathbf{c}_r \right)_n = \widehat{c}_{n,r}$, and $\mathbf{u}_r$ has entries $\left( \mathbf{u}_r\right)_m = \left(\mathbf{u}(x_m)\right)_r$.

Thus the stochastic collocation solution to \eqref{eq:parameterized-model-problem} under the polynomial ansatz \eqref{eq:u-ansatz} with conditions \eqref{eq:collocation-conditions} is given by the solution to \eqref{eq:decoupled-linear-systems}. Therefore, as is common in the stochastic collocation community, we focus entirely on solving the model problem
\begin{align}\label{eq:collocation-problem}
  \mathbf{A} \mathbf{c} \approx \mathbf{u},
\end{align}
for the coefficient vector $\mathbf{c}$ and given data $\mathbf{u}$ obtained from non-intrusive queries of the mode problem \eqref{eq:parameterized-model-problem}.

The three major approximations for \eqref{eq:collocation-problem} that we consider are
\begin{enumerate}
  \item Regression/regularization: we enforce \eqref{eq:collocation-problem} in the least-squares sense over the nodal array
  \item Interpolation: we enforce exact equality in \eqref{eq:collocation-problem}
  \item Compressive sensing: under the assumption that the exact solution coefficients $c_{k,n}$ are sparse in $n$, we attempt to recover a sparse solution $\mathbf{c}$ to \eqref{eq:collocation-problem}.
\end{enumerate}
Our focus in this paper is on presentation on types of 'unstructured', high-order collocation approximation methods. We will introduce the above approximation methods and present some recently-developed unstructured mesh designs that complement each method. Our discussion revolves on the following considerations:
\begin{enumerate}
  \item For a given mesh, how does the stability of the problem scale with the approximation order $N$ and the sample size $M$?
  \item For a given mesh, how is the accuracy of the reconstruction affected by the approximation order $N$ and the sample size $M$?
  \item Can we find or generate a mesh for which the approximation problem has `nice' stability or accuracy properties?
\end{enumerate}
These concerns undergird most of our future discussion.

%% file: content/regression.tex
\section{Regression}\label{sec:methods-regression}
Least-squares regression is one of the most classical approaches to collocation approximation, with vast literature for recovery with noisy data \cite{DeVore_Noise_I,DeVore_Noise_II,Noise_book}. In UQ applications, the uncertainty of the input parameters can be treated as random variables $Z$. It is therefore reasonable to assume that all the stochasticity/uncertainty in the model is described by $Z$.
Thus, we shall focus on the least-squares regression with noise-free data $\mathbf{u}$. Such UQ-focused least-squares regression (sometimes called point collocation \cite{Point}) has been widely explored \cite{NonI,NonII,Zhou_Narayan_Xu,Cohen,FabioL2}.

For the least-squares approach we consider here, the gPC coefficients $\widehat{c}_n$ are estimated by minimizing the mean square error of the response approximation, i.e., one finds the least-squares solution $u_N$ that satisfies
\begin{align}\label{eq:least_solution}
  u_N := \argmin_{p\in P(\Lambda)}  \sum_{m=1}^M \left(p(z_m)-u(z_m)\right)^2,
\end{align}
where the index set $\Lambda$ may be any general index set, e.g. $P^d_k$ or $T^d_k$. 
For convenience we also introduce the discrete inner product on parameter space
\begin{align}\label{eq:least_norm}
  \left\langle u, v\right\rangle_M=\frac{1}{M} \sum_{m=1}^M u(z_m)v(z_m)
\end{align}
and the corresponding discrete norm $\|u\|_{M}=\left\langle u, u \right\rangle^{1/2}_M$.

The formulation \eqref{eq:least_solution} is equivalent to requiring that the model problem \eqref{eq:collocation-problem} is satisfied in the following algebraic sense
\begin{align}\label{eq:algebraic formulation}
\mathbf{c}=\argmin_{\mathbf{v}\in \mathbb{R}^{N}} ||\mathbf{A}\mathbf{v}-\mathbf{u}||_2.
\end{align}
Alternatively, the solution to the least-squares problem (\ref{eq:algebraic formulation}) can also be computed by solving an $N \times N$ system (the ``normal equations"):
\begin{align}\label{eq:normal_equation}
\mathbf{\hat{A}} \mathbf{c} = \mathbf{\hat{u}}
\end{align}
with
\begin{subequations}\label{eq:components}
\begin{align}
  \mathbf{\hat{A}} &:= \mathbf{A}^\top\mathbf{A}= \left(M\left\langle \phi_\alpha, \phi_\beta \right\rangle_M\right)_{\alpha,\beta \in \Lambda} \\
  \mathbf{\hat{u}} &:= \mathbf{A}^\top\mathbf{u}= \left(M\left\langle u, \phi_\alpha\right\rangle \right)_{\alpha \in \Lambda}.
\end{align}
\end{subequations}
We will describe three kinds of such unstructured collocation grids, for which the corresponding theoretical analysis has been addressed recently and is under active development.

\subsection{Monte Carlo sampling}\label{sec:methods-regression-mc}
Monte Carlo (MC) sampling is a natural choice for least-squares regression. For example, one generates independent and identically-distributed (iid) samples from a random variable with density $\rho$ (recall that $Z$ is a random variable with density $\rho$), and these samples form the nodal array $\{z_m\}$. This choice of sampling is certainly justifiable: It is straightforward to establish that the discrete formulation \eqref{eq:algebraic formulation} converges to the $L^2_\rho$-optimal continuous formulation as $M\rightarrow \infty$:
\begin{align*}
  \lim_{M \rightarrow \infty} \argmin_{p \in P(\Lambda)} \sum_{m=1}^M \left( p(z_m) - u(z_m)\right)^2 = \argmin_{p \in P(\Lambda)} \E \left( p(Z) - u(Z) \right)^2
\end{align*}
It is thus not surprising that this least-squares formulation with random samples is popular \cite{Point,NonI,NonII,BS_JCP_LS}. In practical computations, the number of design samples $M$ drawn from the input distribution scales linearly with the dimension of the approximating polynomial space $N$; taking $M \simeq c N$ with $c$ between 2 and 3 is a common choice. Thus we desire stability and convergence results under the assumption that the sample set size $M$ scales \textit{linearly} with the approximation space dimension $N$. Such results are not yet definitively available, but below we discuss what has been accomplished in this direction.

Perhaps the simplest example is MC with $Z$ a uniformly-distributed random variable on a compact interval: Consider $\rho$ the uniform measure on $[-1,1]$. The analysis in \cite{Cohen} shows that if one generates $M \sim N^2$ iid MC design samples drawn from $\rho$, then the spectral properties of the least-squares design matrix are controlled, implying stability for recovery with regression.
\begin{theorem}\cite{Cohen}
  Let $\rho$ be the uniform measure on $[-1,1]$, and $\Lambda = \Lambda^T_{1,N-1}$. For any $r > 0$, assume that $\frac{M}{\log M} \geq C r N^2$ for a universal constant $C$. Then
  \begin{align*}
    \mathrm{Pr} \left[ \vertiii{ \mathbf{\hat{A}} - \mathbf{I} } \geq \frac{1}{2} \right] \leq 2 M^{-r,}.
  \end{align*}
\end{theorem}

This stability result can also be used to prove near-best approximation properties. These results are extended in \cite{Cohen_Fabio_2014} to multidimensional polynomial spaces with arbitrary lower index sets\footnote{A set $\Lambda$ is a \textit{lower set} in this paper if for any $\alpha \in \Lambda$, all indices below it also lie in $\Lambda$: I.e., $\left\{\beta \in \N_0^d\, | \, \beta \leq \alpha\right\} \subset \Lambda$ holds for any $\alpha \in \Lambda$. Here, the ordering $\leq$ is the partial ordering on $\N^d_0$}. I.e., let $\rho$ be the the uniform measure on $[-1; 1]^d$, and assume quadratic dependence $M\sim N^2$ as in the univariate case. If $\Lambda$ in \eqref{eq:components} is a lower set, then similar quadratic scaling $M \propto N^2$ ensures stability and near best approximation of the method independent of the dimension $d$. 


If instead one considers sampling from the Chebyshev measure
\begin{align}\label{eq:chebyshev-density}
  \rho_c(z) &= \prod_{j=1}^d \left(1 - \left(z^{(j)}\right)^2\right)^{-1/2}, & z &\in [-1,1]^d
\end{align}
with the associated Chebyshev polynomial basis $\phi_\alpha$, then the quadratic dependence can be weakened to $M \sim N^{\frac{\ln 3}{\ln 2}}$, where $\frac{\ln 3}{\ln 2} \approx 1.58.$ Such a technique was applied to a class of elliptic PDE models with stochastic coefficients, and an exponential convergence rate in expectation was established \cite{Cohen_Fabio_2014}.

The analysis for unbounded state spaces $D$ is less straightforward; for these unbounded domains, some of the tools used to establish the results above cannot directly be applied. Nevertheless, for the univariate exponential density function $\rho(z) = \exp(-z^2)$, the authors in \cite{Tang_Zhou_2014} use a mapping technique in conjunction with weighted polynomials to establish stability. With this choice of $\rho$, the orthogonal family of Hermite polynomials $\phi_n(z) = H_n(z)$ would be the gPC basis choice. Consider approximation not with polynomials, but instead with the Hermite \textit{functions}:
\begin{align*}
  \psi_n(z) = \exp(-z^2/2) \phi_n(z) = \exp(-z^2/2) H_n(z).
\end{align*}
The collocation samples $z_m$ are \textit{not} generated with respect to the density $\rho$. Instead, one first generates samples $\xi_m$ of a uniform random variable on the interval $[-1,1]$ and subsequently maps these to the real line via the transformation
\begin{align*}
  z(\xi) &= \frac{L}{2} \log \frac{1 - \xi}{1 + \xi}, & \xi &\in (-1,1),
\end{align*}
where the scalar $L$ is a free parameter. Thus, the least-squares design matrix $\mathbf{\hat{A}}$ from \eqref{eq:normal_equation} is now given by
\begin{align}\label{eq:hermite-ls}
  \mathbf{\hat{A}} &= \left( M \left\langle \psi_i, \psi_j \right\rangle_M \right)_{i,j = 1,\ldots, N}, & \langle u, v \rangle_M &= \frac{1}{M} \sum_{m=1}^M u\left( z \left( \xi_m \right)\right) v\left( z\left(\xi_m\right)\right).
\end{align}
The authors in \cite{Tang_Zhou_2014} prove stability of the formulation \eqref{eq:hermite-ls}, requiring only linear scaling of $M$ with respect to $N$ (modulo logarithmic factors).
\begin{theorem}[\cite{Tang_Zhou_2014}]
  For any $r > 0$, assume $\frac{M}{\log M} \gtrsim r\, N$ and $L \gtrsim \sqrt{N}$. Then the least-squares design matrix from \eqref{eq:hermite-ls} is stable in the sense
  \begin{align*}
    \mathrm{Pr} \left[ \vertiii{ \mathbf{\hat{A}} - \mathbf{I} } \geq \frac{5}{8} \right] \leq 2 M^{-r,},
  \end{align*}
  The constant $C$ is universal.
\end{theorem}

Because the hermite functions $\psi_m$ are weighted versions of the Hermite polynomials $H_m$, one can consider the above a statement about stability of a weighted least-squares approximation problem.

In this subsection, Monte Carlo grids are themselves randomly generated, so most statements about stability and accuracy of these least-squares formulations are probabilistic in nature, e.g., convergence with high probability or convergence of the solution expectation. In the next subsection we consider one type of deterministically-generated mesh.

\subsection{The Weil points}\label{sec:methods-regression-weil}
\begin{figure}
\begin{center}
  \includegraphics[width=\textwidth]{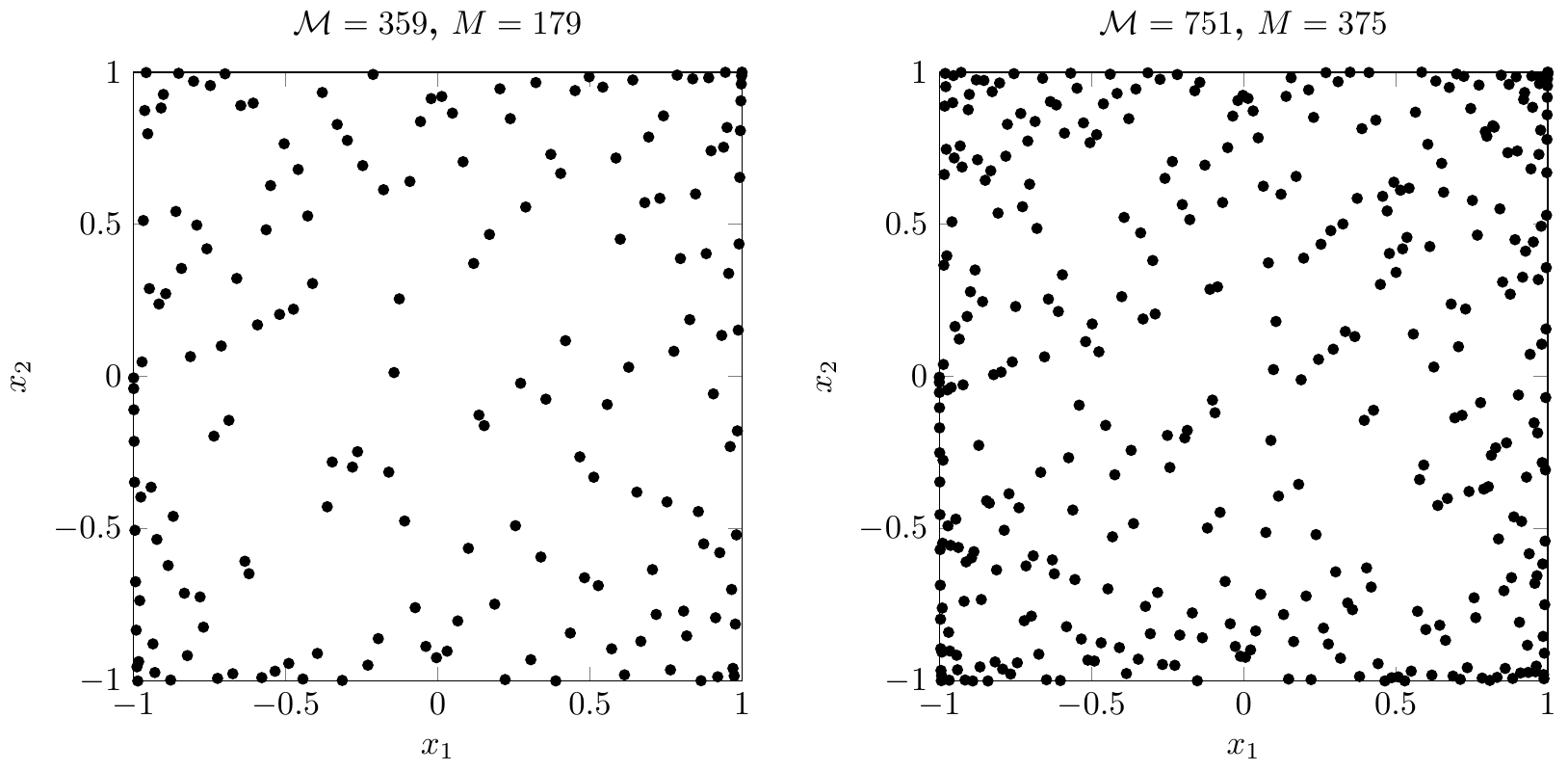}
\end{center}
\caption{Weil grids in two dimensions ($d=2$). Left: prime number seed $\mathcal{M} = 359$ resulting in $M = 179$ points. Right: prime number seed $\mathcal{M} = 751$ resulting in $M = 375$ points.}\label{fig:weil-points}
\end{figure}

In certain applications, a judicious, deterministic choice of samples may provide several advantages over randomly-generated samples. In \cite{Zhou_Narayan_Xu}, the authors present a novel constructive analysis for the deterministically-generated \textit{Weil points}. Suppose that $\mathcal{M}$ is a prime number; \cite{Zhou_Narayan_Xu} proposes the following sample set:
\begin{align}\label{eq:theta-m}
  \mathcal{W}_\mathcal{M}:=\left\{z_{j+1}=\cos( y_{j+1}) : y_j=2\pi\left(j,j^2,\ldots,j^d\right)/\mathcal{M}, \,\,\, j=0,\ldots, \lfloor \mathcal{M}/2 - 1\rfloor \right\},
\end{align}
where $\lfloor \mathcal{M}/2 \rfloor$ gives the integer part of $\mathcal{M}/2.$ Note that the number of the points in the above grid is $M=\lfloor \mathcal{M}/2 \rfloor.$ In fact, with $z_j$ as in \eqref{eq:theta-m}, one can show that the points $\{ z_j\}_0^{M-1}$ coincide with the set of points $\{ z_j\}^\mathcal{M}_{M}$, and so this latter half of the set is discarded. Examples of two-dimensional Weil grids are shown in Figure \ref{fig:weil-points} with $(\mathcal{M}, M) = (359, 179)$ and $(\mathcal{M}, M) = (751, 375)$.

The above collocation grid $\mathcal{W}_M$, designed originally for approximation when $\rho$ is the Chebyshev measure, is motivated by the following formula of Andr\'e Weil (hence the eponymous title ``Weil points"):

\begin{theorem}[Weil's formula \cite{weil_exponential_1948}]\label{thm:weils-formula}
Let $\mathcal{M}$ be a prime number. Suppose $f(w)=m_1w+m_2w^2+\cdots+m_dw^d$ and there is a $j, \,1\leq j\leq d,$ such that $\mathcal{M}\nmid m_j,$ then
\begin{align}\label{eq:summension}
\left| \sum_{j=0}^{\mathcal{M}-1} e^{\frac{2\pi i f(j)}{\mathcal{M}}} \right|\leq (d-1)\sqrt{\mathcal{M}}.
\end{align}
\end{theorem}

This formula plays a central role in deriving several properties about the Weil points, including least-squares stability results. Using Weil's formula, \cite{Zhou_Narayan_Xu} shows that the Weil points distribute asymptotically according to the Chebyshev measure:
\begin{theorem}[\cite{Zhou_Narayan_Xu}]\label{thm:asymptotic-distribution}
  Let $\mathcal{W}_{M_K} = \left\{z_{1,K}, \ldots, z_{M_K, K}\right\}$ be the deterministic sampling set \eqref{eq:theta-m} generated by the $K$'th prime number. 
  We define the empirical measure of the $\mathcal{W}_{M_K}$:
  \begin{align*}
    \nu_K &:= \frac{1}{M_K} \sum_{j=1}^{M_K} \delta_{z_{j,K}}
  \end{align*}
  where $\delta_z$ is the Dirac measure centered at $z$,  and let $\nu_c$ be the Chebyshev measure with density $\dx{\nu_c}(z) = \rho_c(z)$ from \eqref{eq:chebyshev-density}.
  Then $\nu_K \rightarrow \nu_c$ in the weak-$\ast$ toplogy as $K \rightarrow \infty$.
\end{theorem}

Having established what kind of measure the Weil points \eqref{eq:theta-m} sample according to, we turn to stability. Assume that $\rho$ is the Chebyshev density; namely, we use the tensorized Chebyshev polynomials as the gPC polynomial basis. By utilizing the Weil points in the least-squares framework, one can obtain estimates for the components of the design matrix $\mathbf{\hat{A}}$ 
by using Weil's formula Theorem \ref{thm:weils-formula}. These estimates, in conjunction with Gerschgorin's Theorem, result in the following stability result:
\begin{theorem}[\cite{Zhou_Narayan_Xu}]\label{th:stablility}
  Suppose that $\mathbf{I}$ is the size-$N$ identity matrix, $\mathbf{\hat{A}}$ is the design matrix associated with the Chebyshev gPC basis, and the Weil points are generated with the prime number seed $\mathcal{M}$, yielding $M$ points. If $M\geq C(d)\cdot {N}^2,$ then the following stability result holds
\begin{align*}
  \vertiii{\frac{2^{d+1}}{M}\mathbf{\hat{A}}-\mathbf{I}}\leq \frac{1}{2},
\end{align*}
where $|||\cdot|||$ is the spectral norm.
\end{theorem}

Therefore the corresponding least-squares problem admits a unique solution, provided that $M\geq C(d)\cdot {N}^2.$ Stability results for least-squares problems are one ingredient for deriving convergence results. For example, the above stability result yields the following bound on least-squares error:
\begin{theorem}[\cite{Zhou_Narayan_Xu}]\label{th:convergence}
  Let $f\in L^2_\rho$ be a multivariate function, and let $p^\ast_N$ be the $L^\infty$-best approximating polynomial from $P(\Lambda)$, and let $P^N_M f$ be the least-squares $P(\Lambda)$-solution with the Weil points. If the prime number seed satisfies $\mathcal{M}\geq C(d)\cdot {N}^2$, then
\begin{align*}
\|f-P_M^N f\|_{L^2_{\rho}}\leq C \| f- p^\ast_N\|_{L^\infty}.
\end{align*}
\end{theorem}
The above theorem implies the $L^2_\rho$ error for the least-squares solution is comparable to the $L^\infty$ optimal approximation. However, we have assumed the quadratic dependence of the number of design points on the degree of freedom of the approximation space. This is stronger than the dependence $\mathcal{M}\sim N^{\frac{\ln 3}{\ln 2}}$ for the Monte Carlo sampling from the Chebyshev measure as discussed in Section \ref{sec:methods-regression-mc}. Nevertheless the construction here, being deterministic, does not suffer from probabilistic qualifiers on convergence results. (E.g., convergence ``with high probability".) Moreover, we shall show with numerical examples that least-squares regression with Weil points performs comparably to Monte Carlo sampling.

For general (non-Chebyshev) types of orthogonal polynomial approximations on compact domains, \cite{Zhou_Narayan_Xu} also proposed the following weighted least squares approach
\begin{align}\label{eq:weighted_least}
  u_{N,w}= P^N_{M,w} u = \argmin_{p\in P(\Lambda)}  \sum_{m=1}^M w_m\left(u(z_m)-p(z_m)\right)^2,
\end{align}
with positive weights $w_m$ defined by 
\begin{align*}
  w_m = \frac{\rho(z_m)}{\rho_c(z_m)} = \pi^d \rho(z_m) \prod_{q=1}^d \left( 1 - \left(z^{\left(q\right)}_m\right)^2 \right)^{1/2},
\end{align*}
where $z^{\left(q\right)}_m$ is the $q$th component of the grid point $z_m$. The basis set $\phi_\alpha$ is the $\rho$-gPC basis, orthonormal under the weight $\rho$, and $\rho_c$ is the Chebyshev density \eqref{eq:chebyshev-density}. The choice of the weights $w_m$ above ensures that the induced change of measure makes the Weil-sampled $w_m$-weighted norm equivalent to the sought $\rho$-weighted norm.

An example of this will be illustrative: suppose we let $\rho$ be the uniform (probability) density $\rho(z) \equiv 2^{-d}$ on $[-1,1]^d$. Then we have
\begin{align}\label{eq:w-general}
  w_m = \frac{\rho(z_m)}{\rho_c(z_m)} = \left(\pi/2\right)^d \prod_{q=1}^d \left( 1 - \left(z^{(q)}_m\right)^2 \right)^{1/2}.
\end{align}
Since $w_m$ is applied to the quadratic form \eqref{eq:weighted_least}, we are effectively preconditioning $f(z_m)$ with $\sqrt{w_m}$. Thus, if the $\phi_j$ are tensor-product Legendre polynomials (orthonormal under the uniform density), then we are preconditioning our expansion as
\begin{align*}
  \sum_{j=1}^N \widehat{c}_j \phi_j(z) \longrightarrow
  \sum_{j=1}^N \widehat{c}_j \sqrt{w_j} \phi_j(z) =
  \sum_{j=1}^N \widehat{c}_j \left( \prod_{q=1}^d \left(1 - (z^{(q)})^2\right)^{1/4} \phi_j(z)\right).
\end{align*}
This type of preconditioning is known to produce well-conditioned design matrices in the context of $\ell^1$ minimization for Legendre approximations \cite{RW}. Of course if $\rho \propto \rho_c$, then we obtain constant weights. Therefore, the weights proposal \eqref{eq:w-general} reduces to well-known preconditioning techniques for some special cases.


\begin{figure}
\begin{center}
    \includegraphics[width=0.98\textwidth]{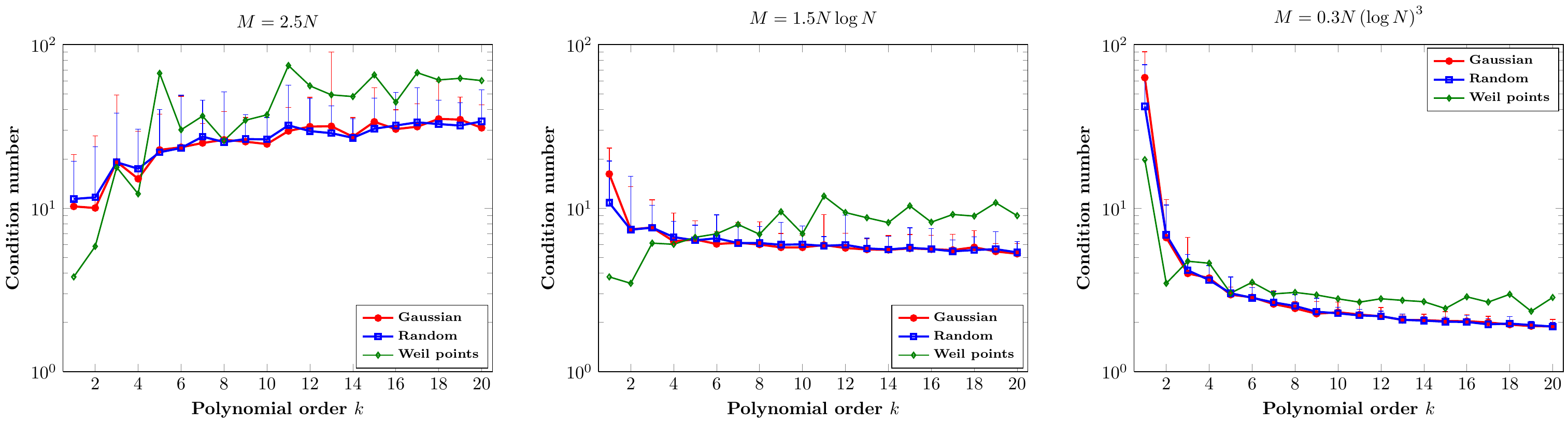}
\end{center}
\caption{Condition number with respect to the polynomial degree $k$ in the 4-dimensional hyperbolic spaces. Left: $M=\frac{5}{2}N.$ Middle: $M=\frac{3}{2}N\log N.$ Right: $M=0.3N\log^3N$}\label{fig:HC4_Condition}
\end{figure}

\subsection{Structured random points}\label{sec:methods-regression-gaussian}
Although a structured grid, such as the tensor-product grid, will frequently be impractical in computations for $d \gg 1$, the grid itself may be used as a candidate set on which to extract a subset of points that is useful for approximation. One way to do this is to randomly sample a subset grid from a high-cardinality structured candidate set, thus producing a subset grid that is essentially unstructrued.


For example, assume that we try to find an approximation in the hyperbolic cross space $H^d_k$  
with $N=h_k^d$.  We assume that the functions $\phi_\alpha(x)$ that span $H_k^d$ are tensor-product Chebyshev polynomials. Such an expansion can also be viewed as an expansion in the tensor product space $P_k^d,$ i.e.
\begin{equation}
u(z)=\sum_{\alpha \in \Lambda^P_{d,k}}  \widehat{c}_\alpha \phi_\alpha(z), \quad  \widehat{c}_\alpha=0\,\,\, \textmd{if}\,\,\,\alpha \in \Lambda^P_{d,k} \backslash \Lambda^H_{d,k}.
\end{equation}
Now let $\{y_m\}_{m=1}^{p_k^d}$ be the tensorized grid of the one dimensional Chebyshev Gauss quadrature points, where $p_k^d=(k+1)^d.$ With such a tensor grid $\{y_m\},$ one can exactly recover any polynomial in $P_k^d$ by the generalized discrete Fourier transform. The normal equations design matrix is the identity matrix in this case. 
Because $h_k^d \ll t_k^d$ when $d$ is large, we have the freedom to choose a subset of points with cardinality $M$ from the full tensor grid, with $M$ satisfying
\begin{equation}
h^d_k  < M < p^d_k.
\end{equation}
We will select these points $z_m$ randomly with the equal probability law on the candidate set $\left\{y_m\right\}$. This idea is proposed in \cite{Zhou_Xiu_2014}, and it is shown that, using the Chebyshev density $\rho$ and associated gPC basis, the design matrix is stable with probability at least $1-2 M^{1-\mu},$ for any $\mu \geq 2,$ provided that
\begin{equation}
\frac{M}{\log M} > C \mu N.
\end{equation}
The framework in \cite{Zhou_Xiu_2014} applies to approximations with densities $\rho$ more general than Chebyshev measures. In particular, it includes measures on bounded domains such as the uniform measure, and on unbounded measures such as the Gaussian measure.

\subsection{Numerical tests}\label{sec:methods-regression-examples}
We now provide some numerical examples to illustrate the stability and the convergence properties of the least-squares approach with the design points described above. Many more related tests are available in \cite{FabioL2,Zhou_Narayan_Xu,Zhou_Xiu_2014}. We first investigate how the number of collocation points affects the condition number,
$
\textmd{cond}(\mathbf{\hat{A}})=\frac{\sigma_{max}(\mathbf{\hat{A}})}{\sigma_{min}(\mathbf{\hat{A}})}
$, where $\sigma_{max}$ and $\sigma_{min}$ are the maximum and minimum singular values, respectively. These results are shown in Figure \ref{fig:HC4_Condition} for the 4-dimensional hyperbolic space $H^4_k$. The orthogonal polynomial measure $\rho$ is the Chebyshev measure and the basis elements are tensor-product Chebyshev polynomials. We test the three design sampling methods described in the previous sections: ``MC" corresponds to the Monte Carlo design of Section \ref{sec:methods-regression-mc}, ``Weil points" corresponds to the Weil points design of Section \ref{sec:methods-regression-weil}, and ``Gaussian" refers to random subset sampling from a tensor-product Gauss quadrature grid from Section \ref{sec:methods-regression-gaussian}.

Because the design matrices with the structured points and the random points are random matrices, we repeat each of these tests 200 times and report the mean condition number. 
We also plotted the error bars for the two kinds of random samples corresponding to one standard deviation above the mean. The left plot of Figure \ref{fig:HC4_Condition} shows results for linear scaling of $M$ with respect to $N$, i.e. $M \sim \frac{5}{2} N,$ and the middle plot shows log-linear dependence $M \sim \frac{3}{2}N \log N$. The right plot shows and even stronger dependence, $M \sim 0.3 N \log^3N,$  We note that the performance of all three methods is similar, in the sense that the log-linear dependence (middle and right plots) admits stable condition numbers with respect to the polynomial order $k$. In contrast, the linear rule (left plot) admits modest growth of the condition number with respect to the polynomial order $k.$ We also observe that the Weil points perform slightly worse compared to the other two types; however, we emphasize that the Weil sampling method and the corresponding analysis is deterministic. Not shown is quadratic dependence, $M \sim N^2$, because the Figure \ref{fig:HC4_Condition} shows that log-linear dependence is enough to guarantee stability. Thus, there seems to be room for improving the quadratic dependence estimate in \cite{Zhou_Narayan_Xu}.
\begin{figure}
\begin{center}
\includegraphics[width=\textwidth]{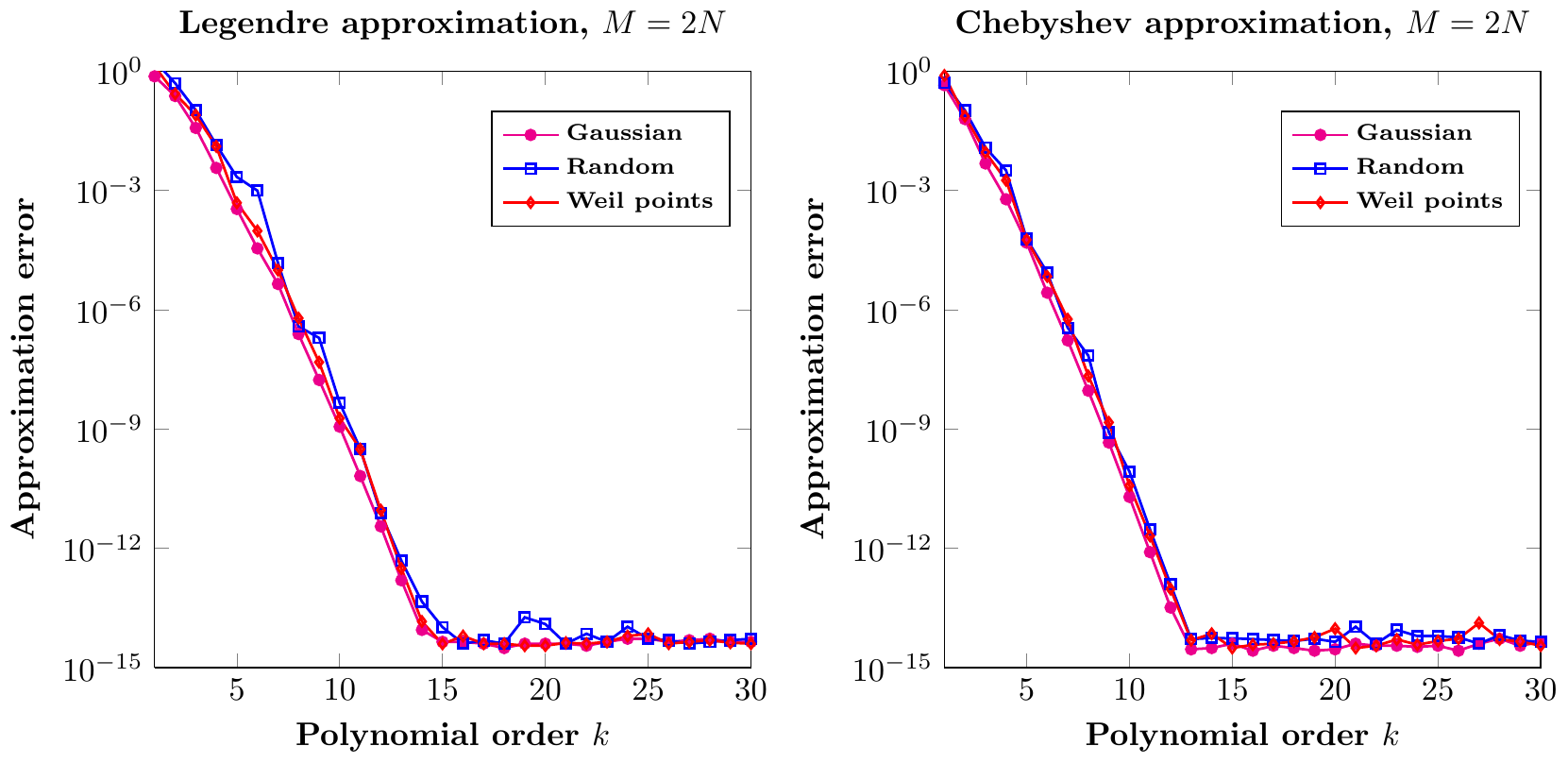}
\end{center}
\caption{Approximation error in the two dimensional total degree space. Left: Chebyshev approximation. Right: Legendre approximation}\label{fig:TD2_Error}
\end{figure}

Finally, we test the convergence rate of the least-squares approaches in $d=2$ dimensions. The target function is smooth, chosen as
$
u(z)=\textmd{exp}^{-\sum_{i=1}^d a_i  z_i},
$
where the parameters $\{a_i\}$ are generated randomly. We measure the error in the $L^\infty$ norm, computed on a set of 2000 points that are independent samples (i.e., different from the design points) from a uniform distribution on $[-1,1]^2.$  The errors with respect to the polynomial order $k$ from the two-dimensional total degree space $T^d_k$ are shown in Figure \ref{fig:TD2_Error}. In the left-hand pane, the underlying measure is the Chebyshev measure and the basis functions $\phi_\alpha$ are Chebyshev polynomials. In the right-hand plot, the underlying measure is the uniform measure and we use the Legendre polynomials as basis elements. In this framework, we will use the preconditioned version of least-squares \eqref{eq:w-general} with Chebyshev-like design points (i.e., random Chebyshev points, the Weil points, and points randomly selected from Chebyshev Gauss points). In all the plots, we have used the linear rule $M = 2 N$ which is dependence that is more feasible in practical computations.

The results given in Figure \ref{fig:TD2_Error} illustrate that the linear rule display the exponential convergence with respect to $k$, for both the least-squares and its preconditioned version. The convergence stagnates at machine precision, which is expected. Such a result also points out a gap between the $M/N$ dependence necessary to achieve
optimality in current theory, and the condition that in practice yields an optimal convergence rate.

\begin{remark}
We note that Quasi Mento Carlo points are deterministically-generated and can also be used in the least-squares framework, and one can find such investigations in \cite{AIAA_QMC,Comparison}. We also remark that the random parameters here are assumed to be supported in $[-1,1]^d.$ One may of course encounter problems with unbounded parameters, e.g., problems with Gaussian/Gamma parameters. Recent work in \cite{Tang_Zhou_2014} investigates such situations.
\end{remark}

%% file: content/cs.tex
\section{Compressive sampling}\label{sec:methods-cs}
Compressive sampling (CS), or compressive sensing, considers recovery the of a sparse representations from limited data, and in applications this is usually considered when there is insufficient information about the target function. In the framework of this paper, this occurs when the number of samples $M$ is less than the cardinality of the polynomial space for the approximation $N$. CS is an emerging and maturing area of research in signal processing that aims at recovering sparse signals accurately from a small number of their random projections (see e.g. \cite{CT_IEEE2005,CT_IEEE2006,CRT_CPAM2006,Donoho2006,Cohen_JAMS} and references therein). A sparse signal is simply a signal that has only few significant coefficients when linearly expanded into a non-adapted basis, in our case the $\{\phi_\alpha\}$. Thus, when the sample count $M$ is smaller than the approximation dimension $N$ then one can use to CS approaches to construct a polynomial approximation, under the assumption that the underlying target function is sparse in $\phi_\alpha$. The typical approach in CS is to minimize a norm of the polynomial, with the constraint of matching the data. For sparsity, one seeks to minimize the $\ell_0$ norm $\left\|\mathbf{c}\right\|_0$ of the coefficient vector, and under certain conditions the minimizing coefficient vector is identical to that which minimizes the $\ell_1$ norm $\left\|\mathbf{c}\right\|_1$. Solving the latter minimization problem is preferred in practice because such problems are convex optimization problems.

The success of the CS lies in the assumption that in practice many target functions are sparse: what appear to be a signal with many features may contain only a small number of notable terms when transformed to the frequency domain. This is indeed the case in many UQ problems. For instance, solutions to linear elliptic PDEs with high-dimensional random coefficients admit sparse representations with respect to the gPC basis under some mild conditions \cite{BS_Sparse,CDS_sparse,TS_Sparse}.

For stochastic collocation methods in the CS framework, we are interested in the case that the number $M$ of solution samples is much smaller than the number $N$ of unknown coefficients. One then seeks to a solution $\mathbf{c}$ with the minimum number of non-zero terms. This can be formulated as the optimization problem
\begin{equation}\label{L0_problem}
\mathbf{P}_0 : \argmin {||\mathbf{c}||_0} \,\,\textrm{subject to}\,\, \mathbf{Ac=u},
\end{equation}
where $||\textbf{c}||_0 :\#\{\alpha: \widehat{c}_\alpha \neq 0\}$ is the $\ell^0$ norm on vectors and should be interpreted as the number of non-zero components of $\mathbf{c}.$ In general, the global minimum solution of $\mathbf{P}_0$ is not unique and is NP-hard to compute.
Fortunately, under \textit{restricted isometry} conditions on the design matrix, the computed $\ell_1$ minimizer approximates the $\ell_0$ minimizer very well, even with noisy measurements \cite{candes_decoding_2005,CT_IEEE2006}. The $\ell_1$ minimization problem in our context is
\begin{equation}\label{eq:L1_problem}
\mathbf{P}_1 : \argmin {||\mathbf{c}||_1} \,\,\textrm{subject to}\,\, \mathbf{Ac=u},
\end{equation}
(Above, $\|\cdot\|_1$ is the $\ell^1$ norm on vectors.) The advantage of the $\mathbf{P}_1$ formulation is that it is a convex problem, and so computational solvers for convex problems may be leveraged \cite{convex_op_I,convex_op_II}.

In practice, since the approximation basis is truncated to a finite number of functions,
the approximation $\mathbf{u}$ may not be exact. Further, the measurement vector $\mathbf{u}$ may be corrupted by noise. These factors lead to the relaxation of the equality in \eqref{eq:L1_problem} and reformulate the problem under the Basis Pursuit Denoising form:
\begin{subequations}
\begin{align}
&\mathbf{P}_{0,\varepsilon} : \argmin{||\mathbf{c}||_0} \,\,\textrm{subject to}\,\,  ||\mathbf{Ac=u}||_2\leq \varepsilon, \label{L0e_problem} \\
\nonumber\\
&\mathbf{P}_{1,\varepsilon} : \argmin{||\mathbf{c}||_1} \,\,\textrm{subject to}\,\, ||\mathbf{Ac=u}||_2\leq \varepsilon,\label{L1e_problem}
\end{align}
\end{subequations}
In what follows, we will focus on the traditional $\ell_1$-minimization form \ref{eq:L1_problem}. Consider the $\ell^p$ error in the best $s$-term approximation of a coefficient vector $\mathbf{c} \in \mathbb{R}^N$
\begin{equation}
\sigma_{s,p}(\mathbf{c}) = \inf_{||\mathbf{y}||_0\leq s} ||\mathbf{y-c}||_p
\end{equation}
Clearly, $\sigma_{s,p}(\mathbf{c}) = 0$ if $\mathbf{c}$ is $s$-sparse, i.e., $||\mathbf{c}||_0 \leq s.$ Let $\mathbf{A}$ be an $M \times N$ matrix. Define the restricted isometry constant (RIC) $\delta_s < 1$ to be the smallest positive number such that the inequality
\begin{equation}
(1-\delta_s)||\mathbf{c}||_2^2 \leq ||\mathbf{A}\mathbf{c}||_2^2 \leq (1+\delta_s)||\mathbf{c}||_2^2
\end{equation}
holds for all $\mathbf{c} \in \mathbb{R}^N$ satisfying $\left\|\mathbf{c}\right\|_0 \leq s$. If the above holds, the matrix $\mathbf{A}$ is said to satisfy the $s$-restricted isometry property (RIP) with restricted isometry constant $\delta_s.$ Now, we are ready to state the following sparse recovery for RIP matrices .
\begin{theorem}[\cite{Rauhut_random_matrix}]\label{thm:cs-rip-error}
Let $\mathbf{A} \in \mathbb{R}^{M\times N}$ be a matrix with RIC such that $\delta_s \leq 0.307.$ For a given $\mathbf{\tilde{c}}\in \mathbb{R}^N,$ let $\mathbf{c}^{\#}$ be the solution of the $\ell_1$-minimization
\begin{equation}\label{eq:direct-minimization}
  \argmin {||\mathbf{c}||_1} \,\,\textmd{subject}\,\, \textmd{to}\,\, \mathbf{Ac=A\tilde{\mathbf{c}}}.
\end{equation}
Then the reconstruction error satisfies
\begin{equation}
||\mathbf{c^{\#}-\tilde{c}}||_2\leq C \frac{\sigma_{s,1}(\tilde{c})}{\sqrt{s}}
\end{equation}
for some constant $C>0$ that depends only on $\delta_s.$ In particular, if $\mathbf{\tilde{c}}$ is $s$-sparse then reconstruction is exact, i.e.,
$\mathbf{c^{\#} = \tilde{c}}.$
\end{theorem}

The theorem above indicates that, as with regression, the choice of nodal array is of great importance: A `good' nodal array will lead to a design matrix $\mathbf{A}$ with an acceptable RIC so that Theorem \ref{thm:cs-rip-error} may be invoked for convergence.

\subsection{Monte Carlo sampling}\label{sec:methods-cs-mc}
As in the least-squares framework, the MC sampling method is still promising in the CS framework. Indeed, much of the pioneering CS work employed MC sampling \cite{CT_IEEE2005,CT_IEEE2006,CRT_CPAM2006,Donoho2006}. However, using such an idea for UQ applications is relatively new.  Some of the first work in this area was investigated in \cite{L1_CICP,L1_JCP} , where the authors applied CS ideas to stochastic collocation and obtained some key properties, such as the probability under which the sparse random response function can be recovered.

The authors in \cite{RW,Rauhut_random_matrix} investigated the recovery of expansions that are sparse in a univariate Legendre polynomial basis. Their analysis relies strongly on RIP results from bounded orthonormal polynomial systems, which we summarize in the following:
\begin{theorem}\cite{Rauhut_random_matrix}\label{th:bounded systems}
Let $\{\phi_n\}$ be a bounded orthonormal system, namely,
\begin{equation}
\sup_n||\phi_n||_\infty = \sup_n\sup_z |\phi_n(z)| \leq L
\end{equation}
for some constant $L\geq 1.$ Let $\mathbf{A} \in \mathbb{R}^{M \times N}$ be the interpolation matrix with entries
$\big\{a_{n,m} = \phi_{n-1}(z_m)\big\}_{1\leq n \leq N, 1\leq m \leq M},$ and let $\mathbf{W}$ be the $M \times M$ diagonal matrix with entries $w_{m,m}=(\pi/2)^{1/2}(1-z^2_m)^{1/4},$  where the points $\{z_m\}_{1\leq m \leq M}$ are iid samples drawn from the one-dimensional Chebyshev measure \eqref{eq:chebyshev-density}. Assume that there is a $\delta > 0$ such that
\begin{equation}
M > C \delta^{-2} L^2 s \log^3(s)\log(N).
\end{equation}
Then with probability at least $1-N^{-\gamma \log^3(s)},$ the RIC $\delta_s$ of  $\frac{1}{\sqrt{M}} \mathbf{W} \mathbf{A}$ satisfies $\delta_s\leq \delta.$ Here the $C, \gamma> 0$ are
universal constants.
\end{theorem}

The authors in \cite{RW} use the above use to provide the following recoverability result for one dimensional sparse Legendre polynomial expansions:
\begin{theorem}[\cite{RW}]\label{thm:cs-rw}
  Let $z_m$ be iid samples from the Chebyshev measure, and let $\mathbf{A} \in \mathbb{R}^{M \times N}$ be the corresponding Legendre polynomial $\phi_n$ design matrix with entries
  \begin{align*}
    \big\{a_{m,n} = \phi_{n-1}(z_m)\big\}_{1\leq n \leq N, 1\leq m \leq M},
  \end{align*}
  and let $\mathbf{W}$ be a diagonal matrix with entries $w_{m,m}=(\pi/2)^{-1/2}(1-z^2_m)^{1/4}.$ Assume that
\begin{equation}
M > C s \log^3(s)\log(N).
\end{equation}
Let $\tilde{\mathbf{c}} \in \R^N$ be any coefficient vector, and consider the solution $\mathbf{c}^{\#}$ to the following $\ell^1$ optimization problem:
\begin{equation}\label{eq:pre-minimization}
\mathbf{c}^{\#}=\argmin {||\mathbf{c}||_1} \quad \textmd{subject}\,\, \textmd{to} \quad \mathbf{WAc=WA\mathbf{\tilde{c}}}.
\end{equation}
Then with high probability
the solution $\mathbf{c}^{\#}$ is within a factor of the best $s$-term approximation error. I.e.,
\begin{align*}
  \textrm{Pr}\left[ \left\| \mathbf{c}^{\#} - \tilde{\mathbf{c}} \right\|_2 \leq \frac{C \sigma_{s,1}\left(\tilde{\mathbf{c}}\right)}{\sqrt{s}} \right] \geq 1-N^{-\gamma \log^3(s)}
\end{align*}
Above, $C$ and $\gamma$ are universal constants.
\end{theorem}

Note that in the above theorem, the random samples are drawn from the Chebyshev measure, namely, the CS method above is a Legendre-preconditioned $\ell_1$-minimization with Chebyshev samples. 
The preconditioned/weighted Legendre polynomials have a uniform bound \cite{Szego} when the weight is $(1 - z^2)^{1/4}$, and this is the main reason to define the weight matrix $\mathbf{W}$ to obtain Theorem \ref{th:bounded systems} or Theorem 4.3.

The above result is extended in \cite{XiuL1} to high-dimensional problems, for both the original $\ell_1$-minimization and the preconditioned $\ell_1$-minimization. We summarize these results with the following theorem:
\begin{theorem}\cite{XiuL1}\label{thm:original_version}
  Let $\{\phi_n\}_{n=0}^{N-1}$ be the Legendre polynomial bases of the total degree space $T_k^d$, and let $u(z)=\sum_{n=0}^{N-1} \widehat{c}_n\phi_n$ be an arbitrary polynomial with coefficient vector $\tilde{\mathbf{c}}$. For some nodal array $\left\{z_m\right\}_{1\leq m \leq M}$, let $\mathbf{A}$ and $\mathbf{W}$ be the corresponding $M \times N$ design matrix and $M \times M$ diagonal preconditioner/weight matrix, respectively, with entries
  \begin{align*}
    a_{m,n} &= \phi_{n-1}(z_m), & w_{m,m} = \left(\frac{\pi}{2}\right)^{-d/2} \prod_{j=1}^d \left( 1 - \left(z^{(j)}_{m}\right)^2 \right)^{1/4}.
  \end{align*}
  \begin{enumerate}
    \item Let $\{z_m\}_{1\leq m \leq M}$ be i.i.d. random samples drawn from the uniform measure, and assume that
      \begin{align*}
      M > 3^k s \log^3(s)\log(N).
      \end{align*}
      Then with high probability, the solution $\mathbf{c}^{\#}$ to the direct $\ell^1$ minimization problem \eqref{eq:direct-minimization} is within a factor of the best $s$-term error:
      \begin{align*}
        \textrm{Pr}\left[ \left\|\mathbf{c}^{\#} -\mathbf{\tilde{c}} \right\|_2 \leq \frac{C\sigma_s(\mathbf{\tilde{c}})_1}{\sqrt{s}} \right] \geq 1 - N^{-\gamma \log^3(s)}
      \end{align*}
    \item Let $\{z_m\}_{1\leq m \leq M}$ be i.i.d. random samples drawn from the Chebyshev measure, and assume that
      \begin{align*}
        M > 2^d s \log^3(s)\log(N).
      \end{align*}
      Then with high probability, the solution $\mathbf{c}^{\#}$ to the preconditioned/weighted $\ell^1$ minimization problem \eqref{eq:pre-minimization} is within a factor of the best $s$-term error:
      \begin{align*}
        \textrm{Pr}\left[ \left\|\mathbf{c}^{\#} -\mathbf{\tilde{c}} \right\|_2 \leq \frac{C\sigma_s(\mathbf{\tilde{c}})_1}{\sqrt{s}} \right] \geq 1 - N^{-\gamma \log^3(s)}
      \end{align*}
  \end{enumerate}
  For both of the above cases, the constants $C$ and $\gamma$ are universal.
\end{theorem}


The above results implies that in very high dimensional spaces $(d>k),$ the Chebyshev preconditioned $\ell^1$-minimization may be less efficient than the direct $\ell^1$-minimization, because it may require more sample points when $2^d > 3^k$. Of course when $d=1$ only the $k=0$ trivial case satisfies this inequality, so in one dimension the preconditioned case is more effective \cite{RW}.

We remark that some modified CS methods have also been investigated in the random sampling framework, one can refer to \cite{L1_Weighted,L1_Reweighted} for the weighted (re-weighted) approaches.

\subsection{Deterministic sampling}\label{sec:methods-cs-weil}

Although random sampling methods have been widely used in the CS framework, a judicious, deterministic choice of points may provide several advantages over randomly-generated points. Deterministic CS sampling and recoverability is a well-studied field for recovery of sparse trigonometric polynomials \cite{Iwen_I,Iwen_II,DeVore,xu,BDFKK_deter}.

However, such a framework has not been widely investigated for UQ applications, especially for recoverability of general types of sparse polynomial expansions. In \cite{XUZHOU}, the authors use the Weil points \eqref{eq:theta-m} to recover sparse Chebyshev polynomials. They use the Weil exponential sum formula, Theorem \ref{thm:weils-formula}, to control the incoherence parameter of the design matrix, which in turn can be used to ascertain RIC information.
More precisely, we have
\begin{theorem}[\cite{XUZHOU}]
Let
  $$
  u =\sum_{\alpha \in \Lambda } c_{\alpha}\phi_\alpha
  $$
be an arbitrary Chebyshev polynomial expansion in $P(\Lambda)$ with coefficient vector $\textbf{c}.$   Suppose that $\mathcal{M} \geq C(\Lambda) s^2$ is a prime number, and assume that $\textbf{c}^{\#}$ is given by the $\ell_1$-minimization problem
\eqref{eq:direct-minimization} with the design matrix $\mathbf{A}$ being the evaluations of the Chebyshev bases on the Weil points generated by prime number seed $\mathcal{M}.$ Then 
\begin{equation}
\|\mathbf{c}^\#-\mathbf{c}\|_2\,\,\lesssim\,\,   \frac{\sigma_{s,1}(\tilde{\mathbf{c}})}{\sqrt{s}}.
\end{equation}
\end{theorem}
Thus, $\ell^1$ minimization can recover $s$-sparse Chebyshev polynomials, provided that the number of the Weil samples scales quadratically with the sparsity $s,$ i.e., $M\sim C(\Lambda) s^2.$ The results apply to any high dimensional polynomial spaces with downward closed multi-index sets, such as the $T^d_k, P^d_k$ and $H^d_k.$ However, different spaces result in different constants $C(\Lambda),$ and the estimates for $C(\Lambda)$ in \cite{XUZHOU} are not optimal.
Although the Weil estimates obtained are not optimal in the sense that they require that the number of sample points $M$ scale quadratically on the sparsity $s$. Nevertheless, we will show via numerical tests that the Weil points have a similar recoverability properties when compared with MC samples. 

The results indicating that the Weil points perform comparably to MC sampling (see Section \ref{sec:methods-cs-examples}) are not necessarily surprising: We have shown in Theorem \ref{thm:asymptotic-distribution} that the Weil points distribute asymptotically according to the Chebyshev measure. Thus, it might be expected that the Weil points produce similar performance as MC Chebyshev samples. Moreover, \cite{XUZHOU} also proposed the preconditioned  $\ell^1$-minimization to handle general sparse polynomials with the Weil points. However, such a framework has the similar drawback as we discussed in the last section: In very high dimensional spaces $(d>k)$, the preconditioned $\ell^1$-minimization may be less efficient than direct $\ell^1$-minimization.
\begin{remark}
  Similar to Section \ref{sec:methods-regression-gaussian}, one may use what we refer to as structured random points (samples randomly chosen from a tensor-product Gauss quadrature grid, or some other candidate set) to recover sparse polynomials. The idea is similar to Section \ref{sec:methods-regression-gaussian}, so we will not discuss it in detail; However, we will test the numerical performance of this method in the next section.
\end{remark}

\subsection{Examples}\label{sec:methods-cs-examples}
\begin{figure}
\begin{center}
\includegraphics[width=\textwidth]{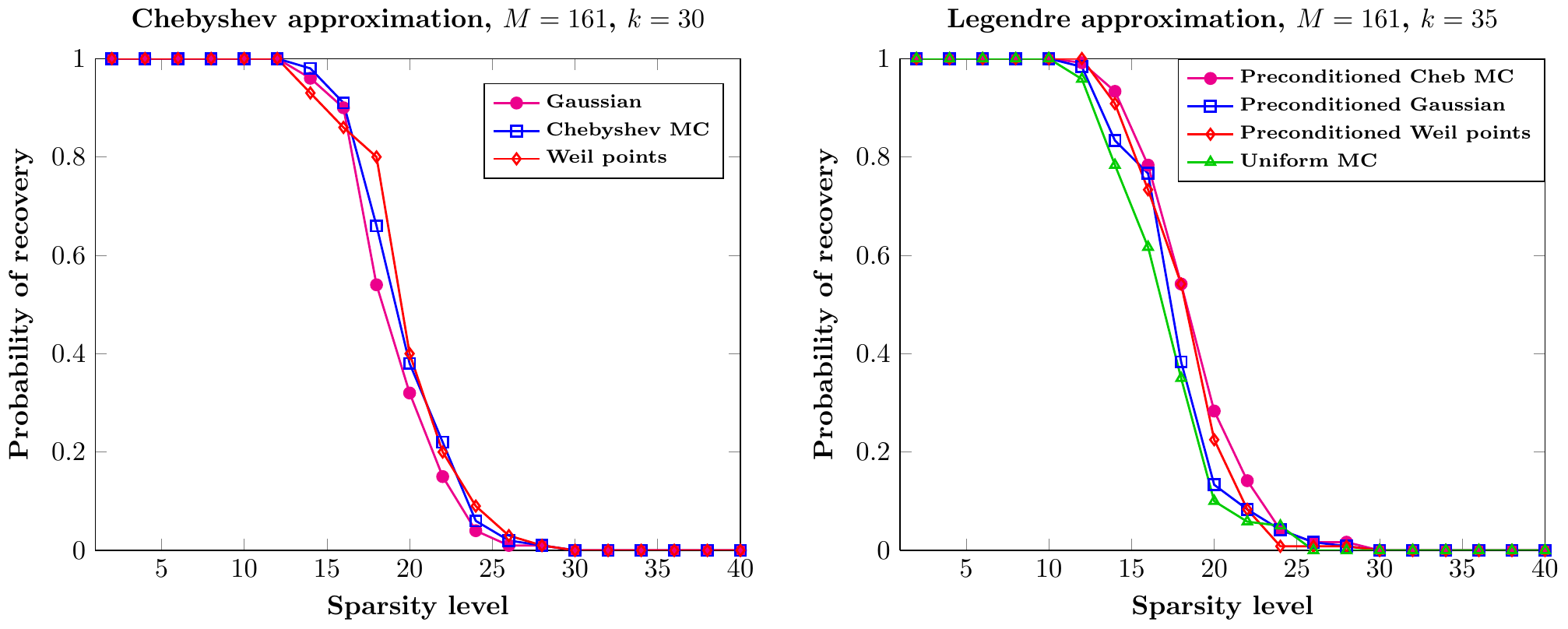}
\end{center}
\caption{Recovery probability with respect to sparsity $s$ in the two-dimensional total degree space $T^2_k$. Left: recovery of sparse Chebyshev polynomial representations. Right: recovery of sparse Legendre polynomial representations.}\label{fig:TD2_recovery}
\end{figure}
This section paralles Section \ref{sec:methods-regression-examples}. We will compare the performance of MC samples, Weil points, and random subsampling from a structured tensor-product Gauss quadrature grid. 
We are interested in the recovery performance with different kinds of points via preconditioned/weighted $\ell_1$ minimization. We assume the target (exact) function has a polynomial form, i.e. $u(z)=\sum_{\alpha \in \Lambda} \widehat{c}_\alpha \phi_\alpha(z)$ with $\left\|\widehat{c}\right\|_0 = s$, and attempt to recover this vector. For a given sparsity level $s$, we fix $s$ coefficients of the polynomial while keeping the rest of the coefficients zero. The values of the $s$ non-zero coefficients are drawn as iid samples from a standard normal distribution. We repeat the experiment $100$ times for each fixed sparsity $s$ and calculate the success rate on these 100 runs. (''Success" here means that $\|\mathbf{\hat{c}-c^{\#}}\|_{\ell^\infty}\leq 10^{-4}.$)

As in Section \ref{sec:methods-regression-examples}, we use the terms ``MC", ``Weil", and ``Gaussian" to refer to Monte Carlo sampling (sections \ref{sec:methods-regression-mc} and \ref{sec:methods-cs-mc}), Weil points sampling (section \ref{sec:methods-regression-weil} and \ref{sec:methods-cs-weil}), and subset sampling from a tensor-product Gauss quadrature grid (section \ref{sec:methods-regression-gaussian}), respectively.

To solve the $\ell_1$ minimization, we employ the available tool SPGL1 from \cite{L1matlab} that was implemented in Matlab. We will conduct two groups of tests, the first using a Chebyshev polynomial expansion, and the second using a Legendre polynomial expansion using the preconditioning strategy of Theorems \ref{thm:cs-rw} and \ref{thm:original_version}.

The first test is a low dimensional test, with the index set $\Lambda$ corresponding to the two dimensional total degree space $T^2_k$. In the left-hand plot of Figure \ref{fig:TD2_recovery}, we show the recovery rate for sparse Chebyshev polynomials (with $k=20, d=2$, and $M=74$) as a function of sparsity level $s$. We see that the three kinds of design points have very similar performance. In the right-hand plot, we show the recovery results for sparse Legendre polynomials (with $k=35, d=2$ and $M=66$) using the preconditioned formulation \eqref{eq:pre-minimization}. We have also tested \textit{direct} (``unpreconditioned") recovery with MC samples drawn from the \textit{uniform} measure. For this low dimensional test, the preconditioned results are slightly better compared to direct $\ell_1$-minimization.

\begin{figure}[h]
\begin{center}
\includegraphics[width=\textwidth]{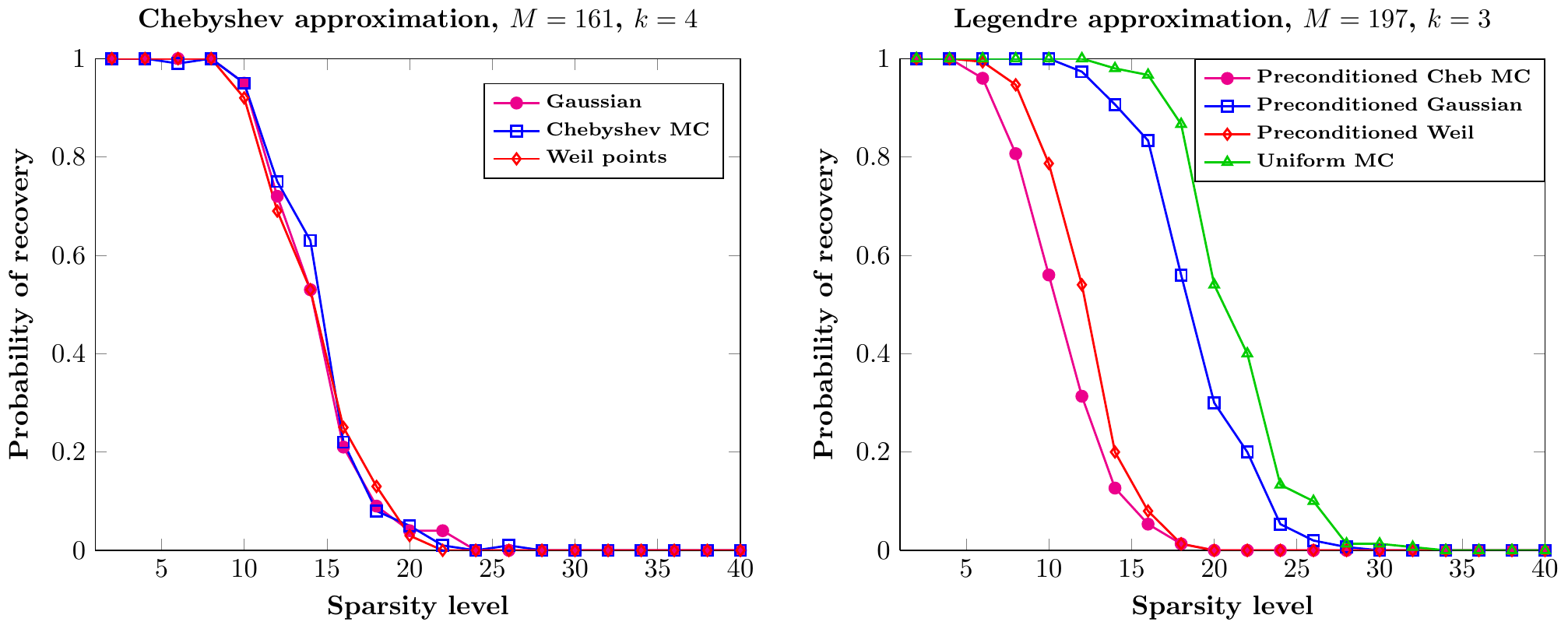}
\end{center}
\caption{Recovery probability with respect to sparsity $s$ in the 15-dimensional total degree space. Left: recovery of sparse Chebyshev polynomial representations. Right: recovery of sparse Legendre polynomial representations.}\label{fig:TD15_recovery}
\end{figure}
The second example we consider is the 15-dimensional total degree space $T^{15}_k$. In Figure \ref{fig:TD15_recovery}, the left-hand plot shows the recovery results for sparse Chebyshev polynomials (with $k=4, d=15$), with respect to the sparsity level $s$ and a fixed number of design points ($M=81$). Again, we see that the three kinds of design points have very similar performance. The right-hand plot provides recovery results for the sparse Legendre polynomials (with $k=4, d=15$, and $M=97$). For this higher dimensional test, the direct $\ell_1$-minimization with MC uniform points performs better than the preconditioned versions (with Weil points, the Chebyshev MC points, and the Gaussian points). This is a direct manifestation of the analysis provided by Theorem \ref{thm:original_version}. Finally, we remark that within the preconditioned tests shown in Figure \ref{fig:TD15_recovery}, the Gaussian points have a noticeably performance than the other two.

%% file: content/interpolation.tex
\section{Interpolation}\label{sec:interpolation}
The last type of approximation we consider is interpolation: given a general array $\mathcal{Z}_M = \left\{z_1, \ldots, z_M\right\}$ of unique nodes and data $u_m = u(z_m)$, we wish to construct a polynomial $p$ satisfying
\begin{align}\label{eq:interpolation-problem}
  p(z) &= \sum_{\alpha \in \Lambda} c_\alpha \phi_\alpha(z), &  p(z_m) &= u_m, & m = 1, \ldots, M,
\end{align}
for some index set $\Lambda$ with $|\Lambda| = M = N$. In the context of the model problem \eqref{eq:collocation-problem}, we have the $N \times N$ system with $\mathbf{A}$ as defined in previous sections:
\begin{align}\label{eq:interpolation-linear-system}
  \mathbf{A} \mathbf{c} &= \mathbf{u}, & \left(A\right)_{n,\alpha} = \phi_\alpha(z_n)
\end{align}
Usually we want $p$ to be as ``simple" as possible, and this translates into the prescription of $\Lambda$, e.g., requiring that $p$ have smallest total degree. One of the challenges for a \naive{} approach to interpolation is that one may not have existence and uniqueness: the classical Mairhuber-Curtis theorem \cite{mairhuber_haars_1956,curtis_parameter_1959} ensures that no matter which $\Lambda$ is chosen, there will exist a set of points on which interpolation is not unisolvent for $P(\Lambda)$.

The above is in contrast to the univariate interpolation problem where, with $N$ nodes, $\deg p$ is known to be exactly $N-1$ regardless of (distinct) nodal choice, and the resulting polynomial is unique given data. In the multivariate case, the space from which the interpolant is drawn cannot be \textit{a priori} specified if we allow the grid $z_m$ to be freely defined; the approximation space must be idenified once the nodes $z_n$ and/or the data $u_n$ is provided.



In Section \ref{sec:interpolation:least-orthogonal} we present one way that uniquely chooses a smallest-degree interpolatory polynomial given data \textit{locations} $z_m$, which depends on the gPC weight function $\rho$ that determines the orthogonal basis $\phi$. The construction is \textit{independent} of the data $\mathbf{u}$. This polynomial is the Least Orthogonal Interpolant and allows us to make smallest-degree polynomial interpolation well-defined in multiple dimensions. This is followed in section \ref{sec:methods-interpolation-lebesgue} by a discussion of Lebesgue constants for accuracy and stability. Finally, sections \ref{sec:interpolation:unweighted} and \ref{sec:interpolation:weighted} describes conditions under which an unstructured nodal array may have well-behaved Lebesgue constant.

\subsection{Least Orthogonal Interpolation}\label{sec:interpolation:least-orthogonal}
Our task in this section is to define and present an operator that prescribes $\rho$-dependent, smallest-degree, unisolvent polynomial interpolation in the multivariate setting. As previously discussed, this is \textit{not} the only choice of interpolation operator, but it has the advantages of being (1) $\rho$-dependent, (2) automatically constructible simply with evaluations of the basis $\phi$, and (3) is compuated via an extension of standard linear algebraic matrix factorizations.

The Least Orthogonal Interpolant (LOI) is a generalization of the `least interpolation' construction \cite{boor_computational_1992,de_boor_least_1992,boor_gauss_1994}. The LOI uses information about the nodal array $z_m$ along with the density $\rho$ to construct this interpolant \cite{narayan_stochastic_2012}. 
If $\rho$ is the density function corresponding to a multivariate standard normal random variable, then the LOI coincides with the traditional least interpolant.

Given a multivariate nodal distribution $z_m$, define $L^2_\rho$ expansions of Dirac distributions centered at the $z_m$:
\begin{align}\label{eq:dirac-measure}
  \delta_{z_m}(\cdot) = \sum_{|\alpha| = 0}^\infty \phi_\alpha\left(z_m\right) \phi_\alpha(\cdot).
\end{align}
Since $\delta_{z_m}$ is a representation of the reproducing kernel, one can directly verify that the formal sum on the right-hand side is the $L^2_\rho$-representor for point-evaluation: $u(z_m) = \left\langle u, \delta_{z_m}\right\rangle_\rho$ for any continuous $u$ in $L^2_\rho$. For any polynomial $p$, one can define the degree-$k$ $L^2_\rho$ projection operators:
\begin{align*}
  \mathcal{P}_k p = \sum_{|\alpha| \leq k} \left\langle p, \phi_\alpha \right\rangle_\rho \phi_\alpha
\end{align*}
Note that these are effectively the action of a series truncation of $\delta_{z_m}$ on $p$. Finally, we introduce the `least-$\rho$' operation:
\begin{align*}
p_{\downarrow,\rho} &= P_{\widehat{k}} p, & \widehat{k} &= \min\left\{k \in \N \; | \; P_k \neq 0\right\}
\end{align*}
So $p_{\downarrow,\rho}$ is effectively the first nonzero ``Taylor series" contribution from an expansion of the form \eqref{eq:dirac-measure}.
Because the projection operators $P_k$ depend on $\rho$, the polynomial $p_{\downarrow,\rho}$ likewise depends on $\rho$. Given a collection of $M$ nodes $\mathcal{Z} = \left\{ z_1, \ldots, z_M \right\}$, the polynomial images of $\delta_{z_m}$ under $(\cdot)_{\downarrow,\rho}$ form a polynomial space:
\begin{align*}
  \Pi_{\mathcal{Z}} = \textrm{span}\left\{ g_{\downarrow,\rho}\, | \, g \in \textrm{span} \left\{ \delta_{z_1}, \ldots, \delta_{z_M} \right\} \right\}.
\end{align*}

The polynomial space $\Pi_{\mathcal{Z}}$ is the least orthogonal polynomial space for interpolation.
\begin{theorem}[\cite{narayan_stochastic_2012}]\label{thm:LOI}
  The space $\Pi_{\mathcal{Z}}$ has dimension $M$, and for any continuous $u$ there is a unique $p \in \Pi_{\mathcal{Z}}$ such that $u_m = p(z_m)$ for $m=1, \ldots, M$. In particular, there exist Lagrange functions $\ell_m(z) \in \Pi_{\mathcal{Z}}$ such that the interpolant $p \in \Pi_{\mathcal{Z}}$ can expressed as
  \begin{align}\label{eq:LOI-lagrange-form}
    p(z) &= \sum_{m=1}^M u_m \ell_m(z), & \ell_m(z_n) &= \delta_{m,n},
  \end{align}
  with $\delta_{m,n}$ the Kronecker delta.
\end{theorem}

The interpolant $p$ is the Least Orthongoal interpolant. It depends on the weight function $\rho$. When $\rho$ is a the Gaussian density function for a standard normal random variable, $\Pi_{\mathcal{Z}}$ coincides with the traditional 'least interpolant' polynomial space of \cite{boor_computational_1992}. While the definition of the space $\Pi_{\mathcal{Z}}$ is fairly abstract, the construction is accomplished in familiar ways: assuming we know $k = \deg \Pi_{\mathcal{Z}}$, let the index set $\Lambda = \mathcal{I}^T_{d,k}$. Form the rectangular design matrix $\mathbf{A}$ from \eqref{eq:decoupled-linear-systems}. The space $\Pi_{\mathcal{Z}}$ can be computed with a combination of $L U$ and $Q R$ operations on this matrix, and the operation count is asymptotically similar to standard interpolatory matrix factorization schemes. Using these operations one obtains the following factorization for the input rectangular matrix $\mathbf{A}$
\begin{align}\label{eq:luh-factorization}
  \bf{P} \bf{A} = \bf{L} \bf{U} \bf{H},
\end{align}
For a size-$M$ interpolation problem, the matrices $\mathbf{L}$ and $\mathbf{U}$ are standard $M \times M$ lower- and upper-triangular matrices corresponding to an $L U$ factorization, and $\mathbf{P}$ is the associated $M \times M$ row permutation matrix. The rectangular matrix $\mathbf{H}$ contains entries that identify the space $\Pi_Z$. 

Once data $\mathbf{u}$ is given, the coefficients $\mathbf{c}$ from \eqref{eq:interpolation-problem} that define the interpolant are computed as
\begin{align*}
  \mathbf{c} = \widetilde{\mathbf{H}}^T \mathbf{U}^{-1} \mathbf{L}^{-1} \mathbf{P} \mathbf{u}
\end{align*}
The new matrix $\widetilde{\mathbf{H}}$ is a rectangular matrix with orthonormal rows that can be obtained from $\mathbf{H}$ in a trivial manner. In order to compute $\mathbf{c}$, the operation count is asymptotically identical to any standard $L U$ factorization for square matrix inversion. We refer to \cite{narayan_stochastic_2012} for the details.


\subsection{Accuracy and the Lebesgue constant}\label{sec:methods-interpolation-lebesgue}
The existence of Lagrange interpolants from Theorem \ref{thm:LOI} allow us to extend univariate notions of stability and accuracy to multivariate settings. With a grid $\mathcal{Z}$ and the associated Lagrange polynomials $\ell_m$ defined in \eqref{eq:LOI-lagrange-form}, we can define a $\mathcal{Z}$-dependent, $\rho$-weighted Lebesgue constant:
\begin{align}\label{eq:weighted-lebesgue-constant}
  \Lambda_\rho(Z) = \max_{z \in \Omega} \rho(z) \sum_{m=1}^M \frac{\left| \ell_m(z) \right|}{\rho(z_m)}.
\end{align}
In the unweighted case $\rho \equiv 1$, this reverts to the more familiar formula for the Lebesgue constant. The introduction of the weight $\rho$ in the above expression is the proper way to generalize the Lebesgue constant to the weighted interpolation problem \cite{lubinsky_survey_2007}. Note that one need not use the LOI space $\Pi_{\mathcal{Z}}$ in the above definition of $\Lambda_\rho$.

On compact domains, the presence of a weight $\rho$ bounded away from $0$ and $\infty$ does not significantly affect the interpolation problem in the context of \eqref{eq:weighted-lebesgue-constant}, as in such a case one can declare the weighted norm equivalent to the unweighted one with appropriate bounding constants. However, on unbounded domains it becomes a necessary consideration. In particular we will consider exponential weights of the form $\rho(z) = \exp( - |z|^\nu)$ for $\nu \geq 1$ on $D = \R^d$.

The error in an interpolation procedure can be understood in terms of the Lebesgue constant, which is the operator norm for weighted interpolation. Let $\mathcal{I}_N$ be any interpolation operator; the LOI construction from the previous section is one such choice. We assume that $\mathcal{I}_N$ maps continous functions to some polynomial space $P_N$, i.e., $\mathcal{I}_N: C_\rho(D) \rightarrow P_N$. The classical Lebesgue Lemma ties the operator norm $\Lambda_\rho$ of $\mathcal{I}_N$ to the interpolant error:
\begin{align}\label{eq:lebesgue-lemma}
  \sup_{z \in D} \rho(z) \left| u(z) - \mathcal{I}_N u(z) \right| \leq \left( 1 + \Lambda_\rho \right) E_{\rho,N}(u).
\end{align}
The term $E_{\rho,N}(u)$ is the best (smallest) possible error in the $\rho$-weighted supremum norm when approximating $u$ by any element in $P_N$. The Lebesgue Lemma is one way of ``separating" the total interpolation error into two parts: first the inherent error $E_{\rho,N}$ depending on the data $u$ and the choice of approximation space, and second the instability $\Lambda_{\rho}$ due to the interpolation procedure. Thus, we are ultimately interested in constructing nodal arrays $\mathcal{Z}$ with small Lebesgue constant. In the multivariate setting, this is still an open problem. 

Values for acceptable Lebesgue constants in one dimension are well-understood: on bounded intervals $D \subset \R$ it is well-known that if a nodal set $\mathcal{Z}$ is uniformly distributed on the interval, then $\Lambda$ grows exponentially as the cardinality of $\mathcal{A}$ is increased (see \cite{trefethen_two_1991} and references therein), and thus the interpolation problem is unstable. On the other hand, for certain special configurations of $\mathcal{A}$ that cluster nodes towards the boundary like the arcsine density $\rho_c = (1 - z^2)^{-1/2}$, the Lebesgue constant exhibits logarithmic growth in the mesh cardinality, yielding a stable interpolation scheme \cite{brutman_lebesgue_1978,gunttner_evaluation_1980}. Results on the one-dimensional unbounded real line are likewise well-established \cite{matjila_bounds_1994,matjila_bounds_1995,szabados_weighted_1997}.

In our multivariate setting, we will generally be concerned with whether or not an array of nodes exhibits exponentially-growing Lebesgue constant. In fact, since $E_{\rho,N}(u)$ usually decreases exponentially with $N$ for $u$ an analytic function, then \textit{subexponential} growth of $\Lambda_\rho$ can be used in Lebesgue's Lemma \eqref{eq:lebesgue-lemma} to conclude convergence of the interpolation scheme.

While concrete results about $\Lambda_\rho$ are available for the univariate case, these results are currently lacking for the multivariate case. What is known are necessary conditions that a grid should satisfy in order to have a ``good" Lebesgue constant: The following sections introduce these conditions, which stem from pluripotential theory. All these results do not directly appeal to the geometry of the grid, and so they apply to any unstructured grids.

\subsection{Unweighted polynomial approximation}\label{sec:interpolation:unweighted}
There is a limited understanding of how one constructs robust interpolation grids in high dimensions. For canonical domains, many approaches yield effective sampling meshes (e.g., the `Padua points' in two dimensions \cite{caliari_bivariate_2005,bos_bivariate_2006}). Tensor-product domains yield good results when the dimension is small ($d \leq 3$). Mapping techniques can produce good nodes for curvilinear domains and blending techniques work well for triangular and tetrahedral domains \cite{warburton_explicit_2006}. Approaches that fit more within the unstructured framework are directly optimized Fekete nodes \cite{taylor_algorithm_2000,bendito_estimation_2007}, approximate Fekete points \cite{sommariva_computing_2009,bos_calculation_2008}, and discrete Leja sequences \cite{bos_computing_2010}. The interested reader may consult these references to find discussions of how to construct multivariate interpolation grids; here we will only discuss theoretical ways of characterizing a `good' interpolation grid.

Consider the case where $D$ is a compact domain in $\R^d$ and $\rho$ is the uniform weight function. Thus, this reduces to unweighted interpolation where the weight appearing in \eqref{eq:weighted-lebesgue-constant} may be set to unity. We are interested in characterizing limiting behaviors of nodal sets $\mathcal{A}$ as $N = |\mathcal{Z}| \rightarrow \infty$. Let $\mathcal{Z}_N$ denote a nodal set of size $N$. We do not require monotonicity, i.e. we do not enforce $\mathcal{Z}_{N} \subset \mathcal{Z}_{N+1}$ for any $N$. We recall that $t^d_k$ is the dimension of the total-degree polynomial space from \eqref{eq:space-dimensions}. The following discussion requires that we restrict our attention to nodal sets with total-degree cardinality $N = t^d_k$ for some degree $k$. We are concerned with the asymptotic behavior of selecting a `good' nodal set $Z_N$; i.e., we are concerned with the behavior of stability metrics of $\mathcal{Z}_N$ as $N \rightarrow \infty$. 


The following notation will be needed: $\mathcal{Z}_N = \left\{z_{1,N}, \ldots, z_{N,N}\right\}$ is a set of $N$ nodes in $D$. We define the following constant as the sum the polynomial degrees of any basis for $T^d_k$:
\begin{align}\label{eq:s-definition}
  s^d_k &\triangleq \sum_{j=1}^k j\, t^{d-1}_{j}, & t^0_j &\triangleq 1
\end{align}
In this section, we will fix $d \geq 1$, and relate the cardinality $N$ of each grid to the polynomial degree $k$:
\begin{align}\label{eq:Nk}
  N_k \triangleq t^d_k.
\end{align}
Given a degree $k$, a set of nodes $\mathcal{Z}^\ast_{N_k}$ is called a set of \textit{Fekete} nodes for the space $T^d_k$ if it maximizes the design matrix determinant over all possible configurations:
\begin{align*}
  \mathcal{Z}^\ast_{N_k} &\triangleq \argmax_{\mathcal{Z} = (y_1, \ldots, y_{N_k}) \subset D} \left| \det \mathbf{A}\left(\mathcal{Z}\right) \right|.
\end{align*}
Above, $\mathbf{A}$ is the (square) design matrix with the \textit{monomial} basis $z^\alpha \in T^d_k$ on the nodes $y_n$. The behavior of the maximal determinant achieved by $\mathcal{Z}^\ast_{N_k}$ defines the transfinite diameter:
\begin{align*}
  \delta(D) \triangleq \lim_{k \rightarrow \infty} \left| \det \mathbf{A}\left(\mathcal{Z}^\ast_{N_k}\right) \right|^{1/s^d_k}.
\end{align*}
For compact $D \subset \R^d$, this limit exists. An array of nodes $\left\{\mathcal{Z}_N\right\}_N$ whose determinant (raised to the $1/s^d_k$) limits to $\delta$ is called \textit{asymptotically} Fekete. 

One final concept we will need is that of equilibrium measures. For a compact $D \subset \R^d$, we introduce the pluripotential equilibrium measure $\mu_{D}$ of the set $D$ \cite{saff_logarithmic_1997,klimeck_pluripotential_1991}. This measure is the multidimensional analogue of the univariate extremal logarithmic energy measure. On a (tensor-product) interval, the measure $\mu_D$ is the (tensor-product) arcsine measure.

The Lebesgue constant, the concept of Fekte nodes, and the equilibrium measure are intimately connected. As described in \cite{bloom_polynomial_1992}, the following theorem considers three ways in which one can characterize the behavior of the array of nodes $\mathcal{Z}_N$ for $N = t^d_k$:
\begin{theorem}[\cite{bos_near_1981, bloom_polynomial_1992, berman_fekete_2011}]\label{thm:unweighted-interpolation}
Consider the three properties an array $Z_{N_k}$ may satisfy:
\begin{enumerate}
  \item \textit{Subexponential Lebesgue constant growth}: $\lim_{k\rightarrow\infty} \left(\Lambda(\mathcal{Z}_{N_k})\right)^{1/k} = 1$
  \item \textit{Asymptotically Fekete}: $\lim_{k\rightarrow\infty} \left| \det \mathbf{A}\left(\mathcal{Z}_{N_k}\right) \right|^{1/s^d_k} = \delta(D)$
  \item \textit{Distribution according to the equilibrium measure}: $$\lim_{k\rightarrow \infty} \frac{1}{N_k} \sum_{n=1}^{N_k} \delta_{z_{n, N_k}} = \mu_{D}$$
\end{enumerate}
Then $1 \Rightarrow 2$ and $2 \Rightarrow 3$; and the reverse implications are false.
\end{theorem}

The proof that $1 \Rightarrow 2$ is given in \cite{bos_near_1981} for one dimension and extends to multiple dimensions as shown in \cite{bloom_polynomial_1992}. That $2 \Rightarrow 3$ is proven for one dimension in \cite{bloom_polynomial_1992} and in \cite{berman_fekete_2011} for the multivariate case. The explanation in \cite{berman_fekete_2011} is relativey abstract, and an accompanying informal discussion can be found in \cite{levenberg_weighted_2010,levenberg_ten_2012}.

These properties give insight into how one should sample an interpolation grid in multiple dimensions: In order to have a stable interpolation operator for interpolation in $T^d_k$ as measured by the Lebesgue constant, one \textit{must} asymptotically sample according to the measure $\mu_{D}$. I.e., sampling according to this measure is a necessary (but not sufficient) condition to obtain a well-bounded Lebesgue constant. 

In the absence of constructive ways to achieve provably well-behaved Lebesgue constants in high dimensions on arbitrary compact domains, one usually devises methods that produce asymptotically Fekete nodes, or nodes that distribute according to $\mu_D$. While such methods are still under active development, promising preliminary results are given by constructing ``approximate" Fekete points \cite{sommariva_computing_2009} or ``discrete" Leja sequences \cite{bos_computing_2010}. A slightly different route may be taken by optimizing a so-called Fejer-constraint that is related to the cardinal Lagrange polynomials of the Least Orthogonal Interpolant \cite{narayan_stochastic_2012}. We also note that the the $d$-dimensional Weil points on a hypercube domain $D$ distribute according to $\mu_{D}$; this is statement of Theorem \ref{thm:asymptotic-distribution}; whether or not these points are asymptotically Fekete is unknown, although use of the Weil points for interpolation is limited by the restriction on their cardinality.


\subsection{Weighted polynomial approximation}\label{sec:interpolation:weighted}
In this section we consider \textit{weighted} approximation on the unbounded domain $D = \R^d$ with exponential weights. We seek to present an extension of Theorem \ref{thm:unweighted-interpolation} to the weighted case; in order to do so we will need \textit{contraction} factors to account for the unbounded nature of $D$.

On $D = \R^d$, we consider weights of the form $\rho(z) = \exp\left( \left\|z \right\|^r \right)$ for $r \geq 1$, where $\|z\|$ is the Euclidean $2$-norm of $z \in \R^d$. This family of weights includes the ubiquitous Gaussian density for $r=2$. It is convenient and common to consider the log-weight, $Q(z) \triangleq - \log \rho = \| z \|^r$. As before, we work on approximation for the space $T^d_k$, define $s^d_k$ as in \eqref{eq:s-definition}, and restrict attention to cardinalities $N_k$ in \eqref{eq:Nk}. Given a set of nodes $\mathcal{Z} = \left\{z_1, \ldots z_{N}\right\}$, then for any $c > 0$ we use the notation $c \mathcal{Z} = \left\{ c z_1, \ldots, c z_N\right\}$. 

The weighted Lebesgue constant $\Lambda_\rho$ is as in \eqref{eq:weighted-lebesgue-constant}.
The notions of Fekete arrays and equilibrium measures carry over to the weighted case with some modifications: an array is a weighted Fekete array for $T^d_k$ if it maximizes a weighted determinant:
\begin{align*}
  \mathcal{Z}^\ast_{N_k,Q} &\triangleq \argmax_{\mathcal{Z} = (y_1, \ldots, y_{N_k}) \subset D} \left| \det \mathbf{A}\left(\mathcal{Z}\right) \right| \prod_{m=1}^{N_k} \rho^{k}(y_j),
\end{align*}
where again the design matrix $\mathbf{A}$ has monomial entries: $(A)_{n,\alpha} = z_n^\alpha$. As with the unweighted case, the limit of the maximal determinants exists \cite{levenberg_ten_2012}; 
\begin{align*}
  \delta_{\rho}(D) \triangleq \lim_{k \rightarrow \infty} \left| \det \mathbf{A}\left(\mathcal{Z}^\ast_{N_k}\right) \prod_{m=1}^{N_k} \rho^{k}(y_j) \right|^{1/s^d_k},
\end{align*}
and any array of points whose weighted determinant behaves like $\delta_\rho(D)$ is called \textit{asymptotically} weighted Fekete. For our domain $D = \R^d$ we will make use of the $\rho$-\textit{weighted} equilibrium measure $\mu_{D, Q}$ from weighted pluripotential theory \cite{saff_logarithmic_1997,klimeck_pluripotential_1991,bloom_weighted_2003}. We mention that the support $\supp \mu_{D,Q}$ of the equilibrium measure is a compact set (even though $D = \R^d$ is unbounded).

For a degree $k$, the weighted situation requires a contraction factor $k^{-1/r}$, where $r$ is the exponent in the log-weight $Q = \|z \|^r$. Interpolation grids with subexponentially growing weighted Lebesgue constant are \textit{not} asymptotically weighted Fekete arrays, but contracted versions of them are.
\begin{theorem}[\cite{berman_fekete_2011,narayan_stochastic_2012}]\label{thm:weighted-interpolation}
  Consider the four properties an array $\mathcal{Z}_{N_k}$ may satisfy:
\begin{enumerate}
  \item[1] \textit{Subexponential weighted Lebesgue constant growth}: $\lim_{k\rightarrow\infty} \left(\Lambda_\rho(\mathcal{Z}_{N_k})\right)^{1/k} = 1$
  \item[2] \textit{Contracted asymptotically weighted Fekete}: $$\lim_{k\rightarrow\infty} \left| \det \mathbf{A}\left(k^{-1/r} \mathcal{Z}_{N_k}\right) \prod_{n=1}^{N_k} \rho^{k}(z_n) \right|^{1/s^d_k} = \delta_\rho(D)$$
  \item[2a] \textit{Uniform contracted boundedness}: There is compact set $S \supset \supp \mu_{\Omega,Q}$ such that $k^{-1/r} \mathcal{Z}_{N_k} \subset S$ for all $k \in \N$.
  \item[3] \textit{Distribution according to the weighted equilibrium measure}: $$\lim_{k\rightarrow \infty} \frac{1}{N_k} \sum_{n=1}^{N_k} \delta_{k^{-1/r} z_{n, N_k}} = \mu_{D,Q}$$
\end{enumerate}
Then $1 \Rightarrow 2$ and $(2 + 2a) \Rightarrow 3$; and the reverse implications are false.\footnote{The uniform contracted boundedness condition 2a is a technicality. Although we suspect that this condition is ultimately unnecessary, a direct proof that $2 \Rightarrow 3$ seems to not yet be available.} 
\end{theorem}

The weighted case on the whole space $\R^d$ here reveals an interesting twist: We do not directly want to sample with respect to $\mu_{D, Q}$; among the justifying reasons is that the measure has compact support, which will be of limited use in approximation with polynomials on an unbounded domain. However, if we sample points so that grids that are progressively contracted by $k^{-1/r}$ distribute according to $\mu_{D,Q}$, then this is a necessary condition for controlled Lebesgue constant growth.


To our knowledge it is unknown whether multidimensional \textit{weighted approximate/discrete} Fekete/Leja arrays distribute according to the weighted pluripotential equilibrium measure. However, we anticipate that such a result is true, and if so would provide a powerful computational approach for generating optimal unstructured meshes. The one-dimensional affirmative answer to this question from \cite{narayan_adaptive_2014} gives hope to this possibility.

%
%
%

%% file: content/conclusion.tex
\section{Summary}
The ability to construct polynomial approximations to functions of high-dimensional inputs is of great interest and importance in modern uncertainty quantification settings. We have explored recent advances in non-intrusive stochastic collocation methods for doing this using a selection of geometrically unstructured meshes. With respect to $N$ degrees of freedom in the polynomial expansion, and $M$ data points, we have presented results for least-squares regularization for overdetermined systems ($M > N$), interpolatory reconstruction for square systems ($M=N$), and compressive sampling/sparse reconstruction for underdetermined systems ($M < N$). 

The results in this field as presented here are far from optimal and settled. We expect a great deal of upcoming and future work on the theory of approximation on unstructured meshes will make such an approach one of the preferred methods for collocative approaches in high-dimensional approximation.

%% file: collocation-review.bbl
\begin{thebibliography}{10}

\bibitem{babuska_stochastic_2010}
I.~Babuska, F.~Nobile, and R.~Tempone.
\newblock A stochastic collocation method for elliptic partial differential
  equations with random input data.
\newblock {\em {SIAM} Review}, 52(2):317--355, January 2010.

\bibitem{barthelmann_high_2000}
V.~Barthelmann, E.~Novak, and K.~Ritter.
\newblock High dimensional polynomial interpolation on sparse grids.
\newblock {\em Advances in Computational Mathematics}, 12(4):273--288, March
  2000.

\bibitem{bellman_adaptive_1961}
R.~E. Bellman.
\newblock {\em Adaptive control processes: a guided tour}.
\newblock Princeton University Press, 1961.

\bibitem{bendito_estimation_2007}
E.~Bendito, A.~Carmona, {A.M.} Encinas, and {J.M.} Gesto.
\newblock Estimation of {F}ekete points.
\newblock {\em Journal of Computational Physics}, 225(2):2354--2376, August
  2007.

\bibitem{berman_fekete_2011}
R.~Berman, S.~Boucksom, and D.~Nystr\"om.
\newblock Fekete points and convergence towards equilibrium measures on complex
  manifolds.
\newblock {\em Acta Mathematica}, 207(1):1--27, 2011.

\bibitem{BS_Sparse}
M.~Bieri and C.~Schwab.
\newblock Sparse high order {FEM} for elliptic {sPDEs}.
\newblock {\em Computer Methods in Applied Mechanics and Engineering},
  198(13-14):1149--1170, March 2009.

\bibitem{DeVore_Noise_II}
P.~Binev, A.~Cohen, W.~Dahmen, and R.~DeVore.
\newblock Universal algorithms for learning theory. ii. piecewise polynomial
  functions.
\newblock {\em Constr. Approx.}, 2(26):127--152, 2007.

\bibitem{DeVore_Noise_I}
P.~Binev, A.~Cohen, W.~Dahmen, R.~A. Devore, and V.~N. Temlyakov.
\newblock Universal algorithms for learning theory part i : Piecewise constant
  functions.
\newblock {\em Journal of Machine Learning Research}, 6:1297--1321, 2005.

\bibitem{BS_JCP_LS}
G.~Blatman and B.~Sudret.
\newblock Adaptive sparse polynomial chaos expansion based on least angle
  regression.
\newblock {\em Journal of Computational Physics}, 230(6):2345--2367, March
  2011.

\bibitem{bloom_polynomial_1992}
T.~Bloom, L.~Bos, C.~Christensen, and N.~Levenberg.
\newblock Polynomial interpolation of holomorphic functions in $\mathbbm{C}$
  and $\mathbbm{C}^n$.
\newblock {\em Rocky Mountain Journal of Mathematics}, 22(2):441--470, June
  1992.

\bibitem{bloom_weighted_2003}
T.~Bloom and N.~Levenberg.
\newblock Weighted pluripotential theory in \$c{\textasciicircum}n\$.
\newblock {\em American Journal of Mathematics}, 125(1):57--103, February 2003.

\bibitem{bos_near_1981}
L.~Bos.
\newblock {\em Near optimal location of points for Lagrange interpolation in
  several variables.}
\newblock PhD thesis, University of Toronto, 1981.

\bibitem{bos_bivariate_2006}
L.~Bos, M.~Caliari, S.~De~Marchi, M.~Vianello, and Y.~Xu.
\newblock Bivariate lagrange interpolation at the padua points: The generating
  curve approach.
\newblock {\em Journal of Approximation Theory}, 143(1):15--25, November 2006.

\bibitem{bos_computing_2010}
L.~Bos, S.~De~Marchi, A.~Sommariva, and M.~Vianello.
\newblock Computing multivariate fekete and leja points by numerical linear
  algebra.
\newblock {\em {SIAM} Journal on Numerical Analysis}, 48(5):1984, 2010.

\bibitem{bos_calculation_2008}
L.~P. Bos and N.~Levenberg.
\newblock On the calculation of approximate fekete points: the univariate case.
\newblock {\em Electronic Transactions on Numerical Analysis}, 30:377--397,
  2008.

\bibitem{BDFKK_deter}
J.~Bourgain, S.~Dilworth, K.~Ford, S.~Konyagin, and D.~Kutzarova.
\newblock Explicit constructions of {RIP} matrices and related problems.
\newblock {\em Duke Mathematical Journal}, 159(1):145--185, July 2011.

\bibitem{breidt_measure-theoretic_2011}
J.~Breidt, T.~Butler, and D.~Estep.
\newblock A measure-theoretic computational method for inverse sensitivity
  problems i: Method and analysis.
\newblock {\em {SIAM} Journal on Numerical Analysis}, 49(5):1836--1859, 2011.

\bibitem{brutman_lebesgue_1978}
L.~Brutman.
\newblock On the lebesgue function for polynomial interpolation.
\newblock {\em {SIAM} Journal on Numerical Analysis}, 15(4):694, 1978.

\bibitem{bungartz_sparse_2004}
H.-J. Bungartz and M.~Griebel.
\newblock Sparse grids.
\newblock {\em Acta Numerica}, 13(-1):147--269, 2004.

\bibitem{burkardt_comparison_2009}
J.~Burkardt and M.~Eldred.
\newblock Comparison of non-intrusive polynomial chaos and stochastic
  collocation methods for uncertainty quantification.
\newblock In {\em 47th {AIAA} Aerospace Sciences Meeting including The New
  Horizons Forum and Aerospace Exposition}. American Institute of Aeronautics
  and Astronautics, 2009.

\bibitem{burkardt_pod_2006}
J.~Burkardt, M.~Gunzburger, and H.-C. Lee.
\newblock {POD} and {CVT}-based reduced-order moeling of navier-stokes flows.
\newblock {\em Computer Methods in Applied Mechanics and Engineering},
  196(1-3):337--355, December 2006.

\bibitem{caliari_bivariate_2005}
M.~Caliari, S.~De~Marchi, and M.~Vianello.
\newblock Bivariate polynomial interpolation on the square at new nodal sets.
\newblock {\em Applied Mathematics and Computation}, 165(2):261--274, June
  2005.

\bibitem{CRT_CPAM2006}
E.~Cand\`{e}s, J.~Romberg, and T.~Tao.
\newblock Stable signal recovery from incomplete and inaccurate measurements.
\newblock {\em Comm. Pure Appl. Math.}, 8(59):1207--1223, 2006.

\bibitem{CT_IEEE2005}
E.~Cand\`{e}s and T.~Tao.
\newblock Stable signal recovery from incomplete and inaccurate measurements.
\newblock {\em IEEE Trans. Inform. Theory}, 12(51):4203--4215, 2005.

\bibitem{candes_decoding_2005}
E.J. Candes and T.~Tao.
\newblock Decoding by linear programming.
\newblock {\em {IEEE} Transactions on Information Theory}, 51(12):4203--4215,
  December 2005.

\bibitem{CT_IEEE2006}
E.J. Candes and T.~Tao.
\newblock Near-optimal signal recovery from random projections: Universal
  encoding strategies?
\newblock {\em {IEEE} Transactions on Information Theory}, 52(12):5406--5425,
  December 2006.

\bibitem{Cohen_Fabio_2014}
A.~Chkifa, A.~Cohen, G.~Migliorati, F.~Nobile, and R.~Tempone.
\newblock Discrete least squares polynomial approximation with random
  evaluations-application to parametric and stochastic elliptic pdes.
\newblock {\em EPFL, MATHICSE Technical Report}, 35/2013.

\bibitem{AIAA_QMC}
S.-K. Choi, R.~V. Grandhi, R.~A. Canfield, and C.~L. Pettit.
\newblock Polynomial chaos expansion with {L}atin hypercube sampling for
  estimating response variability.
\newblock {\em {AIAA} Journal}, 42(6):1191--1198, 2004.

\bibitem{Cohen_JAMS}
A.~Cohen, W.~Dahmen, and R.~A. DeVore.
\newblock Compressed sensing and best $k$-term approximation.
\newblock {\em J. Amer. Math. Soc.}, 22(1):211--231, 2009.

\bibitem{Cohen}
A.~Cohen, M.~A. Davenport, and D.~Leviatan.
\newblock On the stability and accuracy of least squares approximations.
\newblock {\em Foundations of Computational Mathematics}, 13(5):819--834, 2013.

\bibitem{CDS_sparse}
A.~Cohen, R.~DeVore, and C.~Schwab.
\newblock Convergence rates of best n-term galerkin approximations for a class
  of elliptic spdes.
\newblock {\em Found. Comp. Math.}, 6(10):615--646, 2010.

\bibitem{curtis_parameter_1959}
P.~Curtis.
\newblock $n$-parameter families and best approximation.
\newblock {\em Pacific Journal of Mathematics}, 9(4):1013--1027, December 1959.

\bibitem{boor_gauss_1994}
C.~de~Boor.
\newblock Gauss elimination by segments and multivariate polynomial
  interpolation.
\newblock In {\em Proceedings of the conference on Approximation and
  computation : a fetschrift in honor of Walter Gautschi}, pages 11--22, West
  Lafayette, Indiana, United States, 1994. Birkhauser Boston Inc.

\bibitem{boor_computational_1992}
C.~de~Boor and A.~Ron.
\newblock Computational aspects of polynomial interpolation in several
  variables.
\newblock {\em Mathematics of Computation}, 58(198):705--727, April 1992.

\bibitem{de_boor_least_1992}
C.~de~Boor and A.~Ron.
\newblock The least solution for the polynomial interpolation problem.
\newblock {\em Mathematische Zeitschrift}, 210(1):347--378, December 1992.

\bibitem{deane_lowdimensional_1991}
A.~E. Deane, I.~G. Kevrekidis, G.~E. Karniadakis, and S.~A. Orszag.
\newblock Low-dimensional mdoels for complex geometry flows: Application to
  grooved channels and circular cylinders.
\newblock {\em Physics of Fluids A: Fluid Dynamics (1989-1993)},
  3(10):2337--2354, October 1991.

\bibitem{debusschere_protein_2003}
B.~J. Debusschere, H.~N. Najm, A.~Matta, O.~M. Knio, R.~G. Ghanem, and O.~P.
  Le~Maitre.
\newblock Protein labeling reactions in electrochemical microchannel flow:
  Numerical simulation and uncertainty propagation.
\newblock {\em Physics of Fluids}, 15:2238, 2003.

\bibitem{DeVore}
R.~A. DeVore.
\newblock Deterministic constructions of compressed sensing matrices.
\newblock {\em Journal of Complexity}, 23(4-6):918--925, August 2007.

\bibitem{Donoho2006}
D.L. Donoho.
\newblock Compressed sensing.
\newblock {\em IEEE Trans. Inform. Theory}, 4(52):1289--1306, 2006.

\bibitem{L1_JCP}
A.~Doostan and H.~Owhadi.
\newblock A non-adapted sparse approximation of pdes with stochastic inputs.
\newblock {\em J. Comput. Phys.}, 8(230):3015--3034, 2011.

\bibitem{ernst_convergence_2012}
O.G. Ernst, A.~Mugler, H.-J. Starkloff, and E.~Ullmann.
\newblock On the convergence of generalized polynomial chaos expansions.
\newblock {\em {ESAIM}: Mathematical Modelling and Numerical Analysis},
  46(02):317--339, 2012.

\bibitem{Comparison}
Z.~Gao and T.~Zhou.
\newblock Choice of nodal sets for least square polynomial chaos method with
  application to uncertainty quantification.
\newblock {\em Communications in Computational Physics}, 16:365--381, 2014.

\bibitem{gerstner_numerical_1998}
T.~Gerstner and M.~Griebel.
\newblock Numerical integration using sparse grids.
\newblock {\em Numerical Algorithms}, 18(3):209--232, January 1998.

\bibitem{Ghanem}
R.~G. Ghanem and P.~D. Spanos.
\newblock {\em Stochastic finite elements: a spectral approach}.
\newblock Springer-Verlag New York, Inc., 1991.

\bibitem{Noise_book}
L.~Gy{\o}rfi, M.~Kohler, A.~Krzy\.{z}ak, and H.~Walk.
\newblock {\em A Distribution-Free Theory of Nonparametric Regression}.
\newblock Springer Series in Statistics, Springer-Verlag, Berlin, 2002.

\bibitem{gunttner_evaluation_1980}
R.~Günttner.
\newblock Evaluation of lebesgue constants.
\newblock {\em {SIAM} Journal on Numerical Analysis}, 17(4):512--520, August
  1980.

\bibitem{Point}
S.~Hosder, R.~W. Walters, and M.~Balch.
\newblock Point-collocation nonintrusive polynomial chaos method for stochastic
  computational fluid dynamics.
\newblock {\em {AIAA} Journal}, 48(12):2721--2730, December 2010.

\bibitem{Iwen_II}
M.~A. Iwen.
\newblock Simple deterministically constructible rip matrices with sublinear
  fourier sampling requirements.
\newblock {\em 43rd Annual Conference on Information Sciences and Systems
  (CISS), Baltimore, MD}, 2009.

\bibitem{Iwen_I}
M.~A. Iwen.
\newblock Combinatorial sublinear-time fourier algorithms.
\newblock {\em Foundations of Computational Mathematics}, 3(10):303--338, 2010.

\bibitem{kennedy_bayesian_2001}
M.~C. Kennedy and A.~O'Hagan.
\newblock Bayesian calibration of computer models.
\newblock {\em Journal of the Royal Statistical Society: Series B (Statistical
  Methodology)}, 63(3):425--464, January 2001.

\bibitem{klimeck_pluripotential_1991}
M.~Klimeck.
\newblock {\em Pluripotential Theory}.
\newblock Oxford University Press, Oxford, 1991.

\bibitem{levenberg_weighted_2010}
N.~Levenberg.
\newblock Weighted pluripotential theory results of berman-boucksom.
\newblock {\em {arXiv:1010.4035}}, October 2010.

\bibitem{levenberg_ten_2012}
N.~Levenberg.
\newblock Ten lectures on weighted pluripotential theory.
\newblock {\em Dolomites Research Notes on Approximation}, 5:1--59, 2012.

\bibitem{lubinsky_survey_2007}
D.~Lubinsky.
\newblock A survey of weighted polynomial approximation with exponential
  weights.
\newblock {\em Surveys in Approximation Theory}, 3:1--105, 2007.

\bibitem{NonII}
{M. Eldred}.
\newblock Recent advances in non-intrusive polynomial chaos and stochastic
  collocation methods for uncertainty analysis and design.
\newblock In {\em 50th {AIAA}/{ASME}/{ASCE}/{AHS}/{ASC} Structures, Structural
  Dynamics, and Materials Conference}, Structures, Structural Dynamics, and
  Materials and Co-located Conferences. American Institute of Aeronautics and
  Astronautics, May 2009.

\bibitem{mairhuber_haars_1956}
J.C. Mairhuber.
\newblock On haar's theorem concerning chebychev approximation problems having
  unique solutions.
\newblock {\em Proceedings of the American Mathematical Society}, 7(4):609,
  August 1956.

\bibitem{marzouk_stochastic_2009}
Y.~Marzouk and D.~Xiu.
\newblock A stochastic collocation approach to bayesian inference in inverse
  problems.
\newblock {\em Communications in Computational Physics}, 6(4):826--847, October
  2009.

\bibitem{L1_CICP}
L.~Mathelin and K.~A. Gallivan.
\newblock A compressed sensing approach for partial differential equations with
  random input data.
\newblock {\em J. Comput. Phys.}, 12(4):919--954, 2012.

\bibitem{matjila_bounds_1994}
D.~M. Matjila.
\newblock Bounds for lebesgue functions for freud weights.
\newblock {\em Journal of Approximation Theory}, 79(3):385--406, December 1994.

\bibitem{matjila_bounds_1995}
D.~M. Matjila.
\newblock Bounds for the weighted lebesgue functions for freud weights on a
  larger interval.
\newblock {\em Journal of Computational and Applied Mathematics},
  65(1-3):293--298, December 1995.

\bibitem{FabioL2}
G.~Migliorati, F.~Nobile, E.~Schwerin, and R.~Tempone.
\newblock Analysis of the discrete $l^2$ projection on polynomial spaces with
  random evaluations.
\newblock {\em Foundations of Computational Mathematics,
  doi:10.1007/s10208-013-9186-4}, 2014.

\bibitem{najm_uncertainty_2009}
H.~Najm.
\newblock Uncertainty quantification and polynomial chaos techniques in
  computational fluid dynamics.
\newblock {\em Annual review of fluid mechanics}, 41(1):35--52, 2009.

\bibitem{narayan_adaptive_2014}
A.~Narayan and J.~Jakeman.
\newblock Adaptive leja sparse grid constructions for stochastic collocation
  and high-dimensional approximation.
\newblock {\em SIAM Journal on Scientific Computing}, 36(6):A2952--A2983, 2014.

\bibitem{narayan_stochastic_2012}
A.~Narayan and D.~Xiu.
\newblock Stochastic collocation methods on unstructured grids in high
  dimensions via interpolation.
\newblock {\em {SIAM} Journal on Scientific Computing}, 34(3):A1729--A1752,
  June 2012.

\bibitem{FabioC3}
F.~Nobile, R.~Tempone, and C.~G. Webster.
\newblock An anisotropic sparse grid stochastic collocation method for partial
  differential equations with random input data.
\newblock {\em {SIAM} Journal on Numerical Analysis}, 46(5):2411--2442, January
  2008.

\bibitem{patera_reduced_2007}
A.~T. Patera and G.~Rozza.
\newblock {\em Reduced Basis Approximation and A Posteriori Error Estimation
  for Parametrized Partial Differential Equations}.
\newblock {MIT}, version 1.0 edition, 2007.

\bibitem{L1_Weighted}
J.~Penga, J.~Hampton, and A.~Doostan.
\newblock A weighted $\ell^1$-minimization approach for sparse polynomial chaos
  expansions.
\newblock {\em http://arxiv.org/abs/1308.0624v1}, 2013.

\bibitem{poette_uncertainty_2009}
G.~Poëtte, B.~Després, and D.~Lucor.
\newblock Uncertainty quantification for systems of conservation laws.
\newblock {\em Journal of Computational Physics}, 228(7):2443--2467, April
  2009.

\bibitem{prudhomme_reliable_2002}
C.~Prudhomme, D.~V. Rovas, K.~Veroy, L.~Machiels, Y.~Maday, A.~T. Patera, and
  G.~Turinici.
\newblock Reliable real-time solution of parametrized partial differential
  equations: Reduced-basis output bound methods.
\newblock {\em Journal of Fluids Engineering}, 124(1):70, 2002.

\bibitem{pulch_generalised_2012}
R.~Pulch and D.~Xiu.
\newblock Generalised polynomial chaos for a class of linear conservation laws.
\newblock {\em Journal of Scientific Computing}, 51(2):293--312, May 2012.

\bibitem{Rauhut_random_matrix}
H.~Rauhut.
\newblock Compressive sensing and structured random matrices.
\newblock {\em Theoretical Foundations and Numerical Methods for Sparse
  Recovery, Fornasier, M. (Ed.) Berlin,New York (DE GRUYTER)}, pages 1--92,
  2010.

\bibitem{RW}
H.~Rauhut and R.~Ward.
\newblock Sparse {L}egendre expansions via $\ell^1$-minimization.
\newblock {\em Journal of Approximation Theory}, 164(5):517--533, May 2012.

\bibitem{NonI}
{S. Hosder}, {R. Walters}, and {M. Balch}.
\newblock Efficient sampling for non-intrusive polynomial chaos applications
  with multiple uncertain input variables.
\newblock In {\em 48th {AIAA}/{ASME}/{ASCE}/{AHS}/{ASC} Structures, Structural
  Dynamics, and Materials Conference}, Structures, Structural Dynamics, and
  Materials and Co-located Conferences. American Institute of Aeronautics and
  Astronautics, April 2007.

\bibitem{saff_logarithmic_1997}
E.~Saff and V.~Totik.
\newblock {\em Logarithmic Potentials with External Fields}.
\newblock Springer, Berlin, 1997.

\bibitem{sommariva_computing_2009}
A.~Sommariva and M.~Vianello.
\newblock Computing approximate fekete points by {QR} factorizations of
  vandermonde matrices.
\newblock {\em Computers \& Mathematics with Applications}, 57(8):1324--1336,
  April 2009.

\bibitem{szabados_weighted_1997}
J~Szabados.
\newblock Weighted lagrange and hermite-fej\'er interpolation on the real line.
\newblock {\em Journal of Inequalities and Applications}, 1997(2):481267, 1997.

\bibitem{Szego}
G.~Szeg\"{o}.
\newblock {\em Orthogonal Polynomials}.
\newblock American Mathematical Society, Providence, RI, 1975.

\bibitem{Tang_Zhou_2014}
T.~Tang and T.~Zhou.
\newblock On discrete least square projection in unbounded domain with random
  evaluations and its application to parametric uncertainty quantification.
\newblock {\em {SIAM} Journal on Scientific Computing}, 36(5):A2272--A2295,
  2014.

\bibitem{tartakovsky_stochastic_2006}
D.~M. Tartakovsky and D.~Xiu.
\newblock Stochastic analysis of transport in tubes with rough walls.
\newblock {\em Journal of Computational Physics}, 217(1):248--259, September
  2006.

\bibitem{taylor_algorithm_2000}
M.~A. Taylor, B.~A. Wingate, and R.~E. Vincent.
\newblock An algorithm for computing fekete points in the triangle.
\newblock {\em {SIAM} Journal on Numerical Analysis}, 38(5):1707--1720, 2000.

\bibitem{TS_Sparse}
R.A. Todor and C.~Schwab.
\newblock Convergence rates for sparse chaos approximations of elliptic
  problems with stochastic coefficients.
\newblock {\em IMA J. Numer. Anal.}, 2(27):232--261, 2007.

\bibitem{trefethen_two_1991}
L.~N Trefethen and J.~A.~C Weideman.
\newblock Two results on polynomial interpolation in equally spaced points.
\newblock {\em Journal of Approximation Theory}, 65(3):247--260, June 1991.

\bibitem{L1matlab}
E.~van~den Berg and M.~Friedlander.
\newblock Spgl1: A solver for large-scale sparse reconstruction.
\newblock {\em http://www.cs.ubc.ca/labs/ scl/spgl1}, 2007.

\bibitem{warburton_explicit_2006}
T.~Warburton.
\newblock An explicit construction of interpolation nodes on the simplex.
\newblock {\em Journal of Engineering Mathematics}, 56(3):247--262, November
  2006.

\bibitem{weil_exponential_1948}
A.~Weil.
\newblock On some exponential sums.
\newblock {\em Proceedings of the National Academy of Sciences of the United
  States of America}, 34(5):204--207, May 1948.
\newblock {PMID:} 16578290 {PMCID:} {PMC1079093}.

\bibitem{Xiu}
D.~Xiu.
\newblock Efficient collocational approach for parametric uncertainty analysis.
\newblock {\em Commun. Comput. Phys}, page 293�09, 2007.

\bibitem{xiu_fast_2009}
D.~Xiu.
\newblock Fast numerical methods for stochastic computations: A review.
\newblock {\em Communications in Computational Physics}, 5(2-4):242--272, 2009.

\bibitem{xiu_numerical_2010}
D.~Xiu.
\newblock {\em Numerical Methods for Stochastic Computations: A Spectral Method
  Approach}.
\newblock Princeton University Press, July 2010.

\bibitem{XiuH}
D.~Xiu and Jan~S. Hesthaven.
\newblock High-order collocation methods for differential equations with random
  inputs.
\newblock {\em {SIAM} Journal on Scientific Computing}, 27(3):1118--1139,
  January 2005.

\bibitem{XiuK2}
D.~Xiu and G.~E. Karniadakis.
\newblock The wiener--askey polynomial chaos for stochastic differential
  equations.
\newblock {\em {SIAM} Journal on Scientific Computing}, 24(2):619--644, January
  2002.

\bibitem{XiuK1}
D.~Xiu and G.~E. Karniadakis.
\newblock Modeling uncertainty in flow simulations via generalized polynomial
  chaos.
\newblock {\em Journal of Computational Physics}, 187(1):137--167, May 2003.

\bibitem{xu}
Z.~Xu.
\newblock Deterministic sampling of sparse trigonometric polynomials.
\newblock {\em Journal of Complexity}, 27(2):133--140, April 2011.

\bibitem{XUZHOU}
Z.~Xu and T.~Zhou.
\newblock On sparse interpolation and the design of deterministic interpolation
  points.
\newblock {\em to appear in {SIAM} Journal on Scientific Computing}, 2014.

\bibitem{XiuL1}
L.~Yan, L.~Guo, and D.~Xiu.
\newblock Stochastic collocation algorithms using $\ell^1$-minimization.
\newblock {\em International Journal for Uncertainty Quantification},
  2(3):279--293, 2012.

\bibitem{convex_op_II}
J.~Yang and Y.~Zhang.
\newblock Alternating direction algorithms for $\ell^1$-problems in compressive
  sensing.
\newblock {\em {SIAM} Journal on Scientific Computing}, 33(1):250--278, January
  2011.

\bibitem{L1_Reweighted}
X.~Yang and G.~E. Karniadakis.
\newblock Reweighted $\ell_1$ minimization method for stochastic elliptic
  differential equations.
\newblock {\em Journal of Computational Physics}, 248:87--108, September 2013.

\bibitem{convex_op_I}
W.~Yin, S.~Osher, D.~Goldfarb, and J.~Darbon.
\newblock Bregman iterative algorithms for $\ell^1$-minimization with
  applications to compressed sensing.
\newblock {\em {SIAM} J. Imaging Sciences}, 1(1):143--168, 2008.

\bibitem{Zhou_Xiu_2014}
T.~Zhou, A.~Narayan, and D.~Xiu.
\newblock Weighted discrete least-squares polynomial approximation using
  randomized quadratures.
\newblock {\em submitted}, 2014.

\bibitem{Zhou_Narayan_Xu}
T.~Zhou, A.~Narayan, and Z.~Xu.
\newblock Multivariate discrete least-squares approximations with a new type of
  collocation grid.
\newblock {\em {SIAM} Journal on Scientific Computing}, 36(5):A2401--A2422,
  January 2014.

\end{thebibliography}
